\documentclass{amsart}

\usepackage{amsmath,amstext,amsthm,amscd,amsopn,verbatim,amssymb}
\usepackage{hyperref}
\usepackage{mathbbol}
\usepackage{graphicx}
\usepackage{color}
\usepackage{pstricks, pstricks-add, pst-node, pst-coil}
\usepackage{multirow}
\usepackage{array}

\newtheorem{thm}{Theorem}
\newtheorem{prop}[thm]{Proposition}
\newtheorem{lemma}[thm]{Lemma}
\theoremstyle{definition}
\newtheorem{defi}[thm]{Definition}
\newtheorem{propdef}[thm]{Proposition/Definition}

\theoremstyle{remark}
\newtheorem{ex}[thm]{Example}
\newtheorem{rem}[thm]{Remark}

\newcommand{\alg}[1]{\mathfrak{{#1}}}
\newcommand{\cn}[1]{\mathbb{C}^{{#1}}}

\newcommand{\co}[2]{\left[{#1},{#2}\right]} % commutator

\newcommand{\eref}[1]{(\ref{#1})} % equation reference
\newcommand{\morphU}{{\mathcal{U} }} 

\newcommand{\p}{\partial}
\newcommand{\Hom}{\mathop{Hom}}

\newcommand{\R}{{\mathbb{R}}}

\DeclareMathOperator{\dv}{div}

\begin{document}
\title{Formality of cyclic cochains}
%\subtitle{with a section by Damien Calaque}

%\author{Damien Calaque}
\author{Thomas Willwacher}

\address{Department of mathematics, ETH Zurich, 8092 Zurich, Switzerland}
\email{thomas.willwacher@math.ethz.ch, damien.calaque@math.ethz.ch}

\contrib[with a Section by]{Damien Calaque}
%\thanks{The research of D.C.~(on leave of absence from Universit'e Lyon 1) is fully supported by the 
%European Union thanks to a Marie Curie Intra-European Fellowship (contract number MEIF-CT-2007-042212). }
 \thanks{This work was partially supported by the Swiss National Science Foundation (grant 200020-105450).}
%\institute{Department of Mathematics, ETH Zurich \\ thomas.willwacher@math.ethz.ch }
\date{}

\begin{abstract}
We prove Kontsevich's cyclic formality conjecture.
\end{abstract}
\maketitle

\section{Introduction}
In his famous paper \cite{kontsevich}, M. Kontsevich has shown the following Theorem (all notions will be defined below).
\begin{thm}[Kontsevich Formality Theorem]
Let $M$ be a smooth manifold. Then there is an $L_\infty$-quasi-isomorphism of differential graded Lie algebras
\[
(T_{poly}^{\bullet+1},0,\co{\cdot}{\cdot}_{S}) \to (D_{poly}^{\bullet+1}(M), d_H, \co{\cdot}{\cdot}_G)
\]
between the polyvector fields on $M$ and the polydifferential Hochschild complex of $C^\infty(M)$.
\end{thm}

Let now $M$ be an oriented $d$-dimensional manifold with volume form $\omega\in \Omega^d(M)$. The volume form defines a degree $-1$ operator $\dv_\omega$ on the Lie algebra of polyvector fields, that is compatible with the Schouten bracket $\co{\cdot}{\cdot}_{S}$.

Furthermore there is a natural action of the cyclic group of order $n+1$ on $D_{poly}^n(M)$, generated by $\Psi\mapsto \sigma\Psi$, such that for any compactly supported functions $a_0,\dots ,a_n$
\[
\int_M a_0(\sigma\Psi)(a_1,\dots,a_n) \omega= (-1)^n\int_M a_1 \Psi(a_2,\dots,a_n,a_0)\omega.
\]
One can see that the subcomplex of invariants $(D_{poly}^\bullet(M))^\sigma$ is closed under the Gerstenhaber bracket and hence also under the Hochschild differential. The cohomology of this subcomplex is called the (polydifferential) cyclic cohomology of $C^\infty(M)$. Kontsevich conjectured the following variant of his Theorem, which will be proven in the present paper.

\begin{thm}[Cyclic Formality Conjecture]
\label{thm:cycform}
Let $M$ be an oriented manifold with volume form $\omega$. There is an $L_\infty$-quasi-isomorphism of differential graded Lie algebras 
\[
(T_{poly}^{\bullet+1}(M)[u], u \dv_\omega, \co{\cdot}{\cdot}_{S}) \to ((D_{poly}^{\bullet+1}(M))^\sigma, d_H, \co{\cdot}{\cdot}_{G}). 
\]
\end{thm}

The important step in the proof is the construction of a local $L_\infty$-quasiisomorphism that behaves well under globalization.

\begin{thm}
\label{thm:main}
For $M=\R^d$ and $\omega$ a constant volume form, there is an $L_\infty$-quasi-isomorphism $\morphU^{cyc}$ solving Theorem \ref{thm:cycform}. 
\end{thm}

In addition, Proposition \ref{prop:properties} will show that it satisfies properties allowing for globalization (essentially properties P1)-P5) of Kontsevich, see \cite{kontsevich}, section 7).

The morphism $\morphU^{cyc}$ will be given as a sum graphs, more precisely Kontsevich graphs, possibly with tadpoles (edges connecting a vertex to itself). To each tadpole edge, one associates a weight one-form $\eta_z$ as defined in eqn. \eref{equ:etadef}, and for each power of the formal variable $u$ one adds one copy of a two-form $\varpi_z$ defined in eqn. \eref{equ:varpidef}. Otherwise everything is as in the Kontsevich case.
 
A different formula has been proposed by Shoikhet \cite{shoikhet}, involving graphs with dashed pairs. We show that our formula agrees with Shoikhet's for divergence free polyvector fields, hence proving Shoikhet's Conjecture 1, see \cite{shoikhet}. %However, in the general case Shoikhet's morphism does not solve Theorem \ref{thm:cycform} (??? check carefully). In particular, his Conjecture 2 is not correct.

\subsection{Structure of the paper}
In section \ref{sec:hochcyc} we recall the basic notions of Hochschild and cyclic cohomology.
Section \ref{sec:linfty} states our conventions regarding $L_\infty$-algebras.
Section \ref{sec:graphs} recalls the definition of Kontsevich's morphism and defines $\morphU^{cyc}$.
In section \ref{sec:cycinv} we essentially show that $\morphU^{cyc}$ is well-defined.
In section \ref{sec:theproof} Theorem \ref{thm:main} is proven using Stokes' Theorem.
In section \ref{sec:cycproof} Theorem \ref{thm:cycform} is derived from this result.
Section \ref{sec:applications} is dedicated to the single application of the above results we know of, namely the classification of closed star products. %, though unfortunately most of them seem rather irrelevant.

\subsection{Acknowledgements}
The results of this paper emerged during a joint effort with Damien Calaque to prove the cyclic formality conjecture. I am very grateful for the many inspiring and fruitful dicussions with him. Damien also essentially found a proof before I came up with the tidier one decribed in this paper. Furthermore he is the author of section \ref{sec:cycproof} about the globalization of the morphism $\morphU^{cyc}$. 

I have invested some effort and time into obtaining the correct signs and prefactors. The paper of Arnal, Manchon and Masmoudi \cite{arnal} helped me a lot in this task. Hence, although I use different (and hopefully simpler) conventions almost everywhere, credit for the signs, if correct, should go to \cite{arnal}.

Finally, I thank my advisor Giovanni Felder for reading the manuscript and valuable comments. 

\section{Hochschild and Cyclic Cohomology}
\label{sec:hochcyc}
The Hochschild cochain complex of an algebra $A$ is, as a graded vector space
\[
C^\bullet(A) = \Hom(A^{\otimes \bullet},A).
\]
The fundamental operation on this space is called ``braces'':
\begin{align*}
C^k(A)\otimes C^{j_1}(A)\otimes\cdots \otimes C^{j_r}(A) &\to C^{k+j_1+\dots+j_r-r}(A) \\
\phi\otimes \psi_1\otimes \cdots\otimes \psi_r &\mapsto \phi\{\psi_1, \dots, \psi_r\} 
\end{align*}
with
\begin{multline*}
\phi\{\psi_1, \dots, \psi_r\}(a_1,\dots,a_l) = \sum_{0\leq i_1 \leq j_1 \leq \cdots \leq i_r\leq j_r \leq l} 
(-1)^{\sum_{k=1}^r i_k(|\psi_k|-1)}
\\
\phi(a_1,\dots,a_{i_1},\psi_1(a_{i_1+1},\dots,a_{j_1}), a_{j_1+1},\dots, \psi_r(a_{i_r+1},\dots,a_{j_r}),\dots,a_l).
\end{multline*}

Using this operation, one can define a Lie bracket (the Gerstenhaber bracket) on $C^{\bullet+1}(A)$ by
\[
\co{\phi}{\psi}_G := \phi\{ \psi \}-(-1)^{kl} \psi\{\phi \}
\]
for $\phi\in C^{k+1}(A)$ and $\psi\in C^{l+1}(A)$.

\begin{rem}
The product $(\phi,\psi)\mapsto \phi\{\psi\}$ is non-associative in general. However, in the special case that $\phi_1\in C^1(A)$, associativity holds:
\[
(\phi_1\{\phi_2\})\{\psi\} = \phi_1\{\phi_2 \{\psi\} \}. 
\]
This gives rise to an action of the Lie algebra $(C^1(A),\co{\cdot}{\cdot}_G)$ on $C^{\bullet}(A)$:
\[
\phi \cdot \psi = \phi\{\psi\}.
\]
This action is \emph{different} from the ``adjoint'' one, given by the Gerstenhaber bracket:
\[
(\phi,\psi)\mapsto \co{\phi}{\psi}_G.
\]
\end{rem}

The multiplication in $A$ defines a canonical element $m\in C^2(A)$. The Hochschild differential is defined as 
\[
d_H \phi = \co{m}{\phi}_G.
\]
The cohomology of $(C^\bullet(A), d_H)$ is called the Hochschild cohomology of $A$. 

There is an associative (cup) product $\cup$ on $C^\bullet(A)$ given by 
\[
\phi\otimes \psi \mapsto \phi\cup \psi := m\{\phi,\psi\}.
\]

\subsection{Polyvector Fields and Polydifferential Operators}
In this paper, the algebra $A=C^\infty(M)$ will be the algebra of smooth functions on a $d$-dimensional manifold $M$. The algebra of polyvector fields on $M$, $T^\bullet_{poly}(M)$ is the algebra of smooth sections of $\wedge^\bullet TM$. 
There is Lie bracket $\co{\cdot}{\cdot}_{S}$ on $T^{\bullet+1}_{poly}(M)$, the Schouten bracket, extending the Lie derivative and making $T^\bullet_{poly}(M)$ a Gerstenhaber algebra. More concretely,
\begin{multline*}
\co{v_1\wedge\cdots \wedge v_m}{w_1\wedge\cdots \wedge w_n}_{S}
=\\= 
\sum_{i=0}^n\sum_{j=0}^n (-1)^{i+j}\co{v_i}{w_j}\wedge v_1\wedge \cdots \wedge\hat{v}_i\wedge\cdots \wedge v_m\wedge w_1\wedge\cdots \wedge \hat{w}_j\wedge\cdots \wedge w_n.
\end{multline*}
For the special case $M=\R^d$ with standard coordinates we also introduce the notation
\[
%\alpha^{i_1\cdots i_n} \p_{i_1}\wedge \cdots \wedge \p_{i_n} \bullet \beta^J \p_J
\gamma_1\bullet \gamma_2 = \sum_{i=1}^d (\iota_{\p_i}\gamma_1)\wedge (\p_i\cdot \gamma_2)
%=
%\sum_{k=1}^n(-1)^{k+1} \alpha^{i_1\cdots i_n} (\p_{i_k}\beta^J) \p_{i_1}\wedge \cdots \wedge\hat{\p}_{i_k}\wedge \cdots \wedge \p_{i_n} \wedge \p_J
\]
where $\iota_{\p_i}$ is the insertion and $\p_i\cdot$ denotes a Lie derivative. With this definition $\co{\gamma_1}{\gamma_2}_{S}= (-1)^{k_1-1}(\gamma_1\bullet\gamma_2+(-1)^{k_1 k_2} \gamma_2\bullet\gamma_1)$ for all $\gamma_1\in T^{k_1}_{poly}(M)$ and $\gamma_2\in T^{k_2}_{poly}(M)$.

Assume now that $M$ is oriented, with volume form $\omega$. Contraction with $\omega$ defines an isomorphism $T^\bullet_{poly}(M)\to \Omega^{d-\bullet}(M)$. The divergence operator $\dv_\omega$ on $T^\bullet_{poly}(M)$ is defined as the pull-back of the de Rham differential $d$ on $\Omega^{\bullet}(M)$ under this isomorphism. One can check that $\dv_\omega$ is a derivation with respect to the Schouten bracket, i.e., %for all $\gamma_1\in T^{k_1}_{poly}(M)$ and $\gamma_2\in T^{k_2}_{poly}(M)$
\[
\dv_\omega \co{\gamma_1}{\gamma_2}_{S} = \co{\dv_\omega\gamma_1}{\gamma_2}_{S}+ (-1)^{k_1-1}\co{\gamma_1}{\dv_\omega \gamma_2}_{S}.
\]
 %We denote the subspace of polyvector fields annihilated by $\dv_\omega$ by $(T^\bullet_{poly}(M))_{div}$.

Let $D_{poly}^1(M)\subset C^1(A)$ be the space of differential operators. We define the space of polydifferential operators $D_{poly}^\bullet(M) \subset C^\bullet(A)$ to be the subcomplex generated by $D_{poly}^1(M)$ and $A=C^0(A)=D_{poly}^0(M)$ under the cup product $\cup$.

The Lie algebra of vector fields $T^1_{poly}(M)$ acts on $D_{poly}^\bullet(M)$ by
\[
(v, \phi) \mapsto v\cdot \phi = v\{\phi\}.
\]
It is easy to check that this action, together with the action by Lie derivatives on $\Omega^d(M)$, gives rise to an action of $T^1_{poly}(M)$ on
\[
\widetilde{CC}^\bullet(M)= D_{poly}^{\bullet+1}(M) \otimes_A \Omega^d(M)
\]
Here the tensor product is defined using the $A$-module structure on $D_{poly}^{\bullet+1}(M)$ given by the cup product.

\subsection{Cyclic Cohomology}
We already defined the action of the cyclic group on an polydifferential operator, involving an integral, functions of compact support and integration by parts. For later proofs it will be convenient to have a different, but equivalent definition.

For this, we replace the complex of polydifferential operators by the coinvariant space 
\[
CC^\bullet(M) = \left( \widetilde{CC}^\bullet(M) \right)_{T^1_{poly}(M)}.
\]
The space $CC^k(M)$ carries a natural action of the cyclic group of order $k+1$, coming from the action on $D_{poly}^{k+1}(M)$ by interchange of the arguments. We call the generator $\sigma$:
\[
(\sigma \phi)(a_0,..,a_k)=(-1)^k\phi(a_1,..,a_k,a_0).
\]
There is a natural map defined as the composition
\[
\iota: D_{poly}^{\bullet}(M) \hookrightarrow D_{poly}^{\bullet+1}(M) 
\stackrel{\cdot \omega}{\to}\widetilde{CC}^\bullet(M) \twoheadrightarrow CC^\bullet(M).
\]
Here the first map is defined as 
\[
\phi \mapsto 1\cup \phi.
\]
where the ``1'' is the differential operator $1\in D_{poly}^{1}(M)$, not to be confused with the function $1\in D_{poly}^{0}(M)=A$.

\begin{lemma}
\label{lem:dccisom}
The map $\iota$ is an isomorphism of graded vector spaces.
\end{lemma}
\begin{proof}
Through ``integration by parts'' each element of $CC^k(M)$ has a representative $A$-linear in its first argument. Those can be identified with the image of $D_{poly}^{k}(M)$ in $\widetilde{CC}^k(M)$, and hence $\iota$ is surjective. Injectivity follows from the fact that if the expression
\[
\sum_j L_{X_j} (f g_j \omega) = \sum_j (f \dv_\omega(g_j X_j)+g_jX_j\cdot f) \omega
\]
contains no derivatives in $f$, then $\sum_j g_jX_j=0$ and it vanishes identically.
\end{proof}

Using this isomorphism, we define the operator $d_H$ and the bracket $\co{\cdot}{\cdot}_G$ also on $CC^\bullet(M)$ and the operator $\sigma$ also on $D_{poly}^\bullet(M)$. It is easily seen that the latter definition of $\sigma$ coincides with the one given in the introduction.

The following Proposition/Definition is essentially due to Connes.
\begin{propdef}
The Hochschild differential $d_H$ leaves invariant the subspace
\[
(D_{poly}^\bullet(M))^\sigma \cong (CC^\bullet(M))^\sigma
\]
of invariants under the cyclic group action. The resulting complex we call the \emph{polydifferential cyclic chain complex} of $A$.
\end{propdef}

As was mentioned in the introduction, this statement is a direct consequence of the more general statement.
\begin{prop}
The space $(D_{poly}^\bullet(M))^\sigma$ is closed under the Gerstenhaber bracket.
\end{prop}

\section{$L_\infty$ algebras}
\label{sec:linfty}
%\subsection{The symmetric coalgebra}
Let $V$ be a graded vector space. We denote the symmetric algebra by $S^\bullet V=\bigoplus_{n\geq 0} V^{\otimes n} / I$, where $I$ is the two-sided ideal generated by relations $x\otimes y -(-1)^{|x||y|}y\otimes x$. The product will be denoted by $\odot$. For example, the expressions $x_1\odot \cdots \odot x_n :=[x_1\otimes \cdots \otimes x_n]$ generate $S^nV$ as a vector space. Let $S^+ V:= \bigoplus_{n\geq 1} S^n V$, with grading given by $|x_1\odot \cdots \odot x_n|=\sum_j|x_j|$. This space carries the structure of a graded cocommutative coalgebra without counit, with comultiplication given by
\[
\Delta( x_1 \odot \cdots \odot x_n) = \sum_{\substack{I\sqcup J=[n] \\ |I|,|J|\geq 1}} \epsilon(I,J) \bigodot_{i\in I} x_i \otimes \bigodot_{j\in J} x_j.
\]
Here $\epsilon(I,J)$ is the sign of the ``shuffle'' permutation bringing the elements of $I$ and $J$ corresponding to \emph{odd} $x$'s into increasing order. Note that $\epsilon(I,J)$ implicitly depends on the degrees of the $x$'s.

\begin{defi}
Let $(\mathcal{C},\Delta)$, $(\mathcal{C}',\Delta')$ be graded coalgebras. A linear map $\mathcal{F}\in \Hom_k(\mathcal{C},\mathcal{C}')$ is called \emph{degree $k$ morphism of coalgebras} if $\Delta' \circ \mathcal{F} = (\mathcal{F}\otimes \mathcal{F})\circ \Delta$. A linear map $Q\in \Hom_k(\mathcal{C},\mathcal{C})$ is called a \emph{degree $k$ coderivation} on $\mathcal{C}$ if  $\Delta \circ Q = (Q\otimes 1 + 1\otimes Q)\circ \Delta$. Here we use the Koszul sign rule, e.g., $(1\otimes Q)(x\otimes y) = (-1)^{k|x|}x\otimes Qy$ etc.
%\[
%\Delta_1 Q x = \sum_{(x)} Qx'\otimes x'' + (-1)^{k|x'|}x'\otimes Qx''
%\]
%using Sweedler notation.
\end{defi}

Any coderivation $Q$ on $S^+V$ (coalgebra morphism $\mathcal{F}:S^+V\to S^+W$) is uniquely determined by its composition with the projection $S^+V\to S^1V=V$ ($S^+W\to S^1W=W$). The restriction to $S^nV$ of this composition will be denoted by $Q_n\in \Hom(S^kV,V)$ ($\mathcal{F}_n\in \Hom(S^kV,W)$) and called the $n$-th ``Taylor coefficient'' of $Q$ ($\mathcal{F}$).

\begin{defi}
An $L_\infty$-algebra structure on a graded vector space $\alg{g}^\bullet$ is a degree 1 coderivation $Q$ on $S^+(\alg{g}^{\bullet+1})$ such that $Q^2=0$. A morphism of $L_\infty$ algebras $\mathcal{F}:(\alg{g},Q)\to (\alg{g}',Q')$ is a degree 0 coalgebra morphism $\mathcal{F}:S^+(\alg{g}^{\bullet+1})\to S^+((\alg{g}')^{\bullet+1})$ such that $\mathcal{F}Q=Q'\mathcal{F}$.
\end{defi}

In components, the $L_\infty$-relations read
\[
\sum_{\substack{I\sqcup J=[n] \\ |I|,|J|\geq 1} } \epsilon(I,J) Q_{|J|+1}(Q_{|I|}(\bigodot_{i\in I} x_i) \odot \bigodot_{j\in J} x_j) = 0.
\]
All $L_\infty$-algebras in this paper will be of the following type:
\begin{ex}
Let $(\alg{g},d,\co{\cdot}{\cdot})$ be a differential graded Lie algebra. Then the assignments $Q_1(x)=d x$, $Q_2(x,y)=(-1)^{|x|}\co{x}{y}$, $Q_n=0$ for $n=3,4,..$ define an $L_\infty$-algebra structure on $\alg{g}$. To see this, calculate
\begin{align*}
& Q_1(Q_2(x\odot y))+Q_2(Q_1(x)\odot y) + (-1)^{|x||y|}Q_2(Q_1(y)\odot x)
\\ &\quad =
(-1)^{|x|}d\co{x}{y} - (-1)^{|x|}\co{dx}{y}+ (-1)^{|x||y|+|y|+1}\co{dy}{x} 
\\ &\quad =
(-1)^{|x|}\left( d\co{x}{y} - \co{dx}{y}- (-1)^{|x|+|x||y|+|y|+1+(|x|+1)|y|}\co{x}{dy} \right)= 0.
\end{align*}
Here and everywhere in the paper $|x|$ is the degree wrt. the grading on the coalgebra, i.e.,  $x\in \alg{g}^{|x|+1}$.
\end{ex}
\begin{ex}
An $L_\infty$-morphism $\mathcal{F}$ between dglas $\alg{g}$, $\alg{g}'$ has to satisfy the relations
\begin{multline}
\label{equ:dglaLinfty}
Q_1'\mathcal{F}_n(x_1\odot \cdots \odot x_n) + \frac{1}{2} \sum_{\substack{I\sqcup J=[n] \\ |I|,|J|\geq 1} }\epsilon(I,J)
Q_2'(\mathcal{F}_{|I|}(\bigodot_{i\in I} x_i) \odot \mathcal{F}_{|J|}( \bigodot_{j\in J} x_j) )
= \\ =
\sum_{i=1}^n \epsilon(i,1,\dots,\hat{i},\dots, n) \mathcal{F}_n( Q_1(x_i)\odot x_1\odot \cdots\odot\hat{x}_i \odot \cdots \odot x_n)
+ \\ +
\frac{1}{2}
\sum_{i\neq j}^n \epsilon (i,j,\dots,\hat{i},\dots,\hat{j},\dots, n) \mathcal{F}_{n-1}( Q_2(x_i\odot x_j)\odot x_1\odot\cdots\odot \hat{x}_i \odot \cdots \odot\hat{x}_j \odot \cdots x_n).
\end{multline}
Here the factor $\epsilon(i,j,1,..,\hat{i},..,\hat{j},.., n)$ is the sign of the permutation on the \emph{odd} $x$'s that brings $x_i$ and $x_j$ to the left.

\end{ex}

\subsection{Special case: $T_{poly}$ and $(D_{poly})^\sigma$}
We consider here the special case $\alg{g}=(T_{poly}(M)[u],u \dv, \co{\cdot}{\cdot}_{S})$, $\alg{g}'=((D_{poly}(M))^\sigma,d_H, \co{\cdot}{\cdot}_{G})$ and $M=\R^d$. 

Let $x_1, \dots , x_n\in T_{poly}(M)[u]$ and denote $|x_j|=k_j$. Then
\begin{align*}
&\frac{1}{2} \sum_{\substack{I\sqcup J=[n] \\ |I|,|J|\geq 1} }\epsilon(I,J)
Q_2'(\mathcal{F}_{|I|}(\bigodot_{i\in I} x_i) \odot \mathcal{F}_{|J|}( \bigodot_{j\in J} x_j) )
\\ &\quad =
\frac{1}{2} \sum_{\substack{I\sqcup J=[n] \\ |I|,|J|\geq 1} }\epsilon(I,J) (-1)^{|k_I|}
\left(
\mathcal{F}_{|I|}(\bigodot_{i\in I} x_i) \circ \mathcal{F}_{|J|}( \bigodot_{j\in J} x_j) 
-(-1)^{(|k_I|+1)(|k_J|+1)} \mathcal{F}_{|J|}( \bigodot_{j\in J} x_j) \circ \mathcal{F}_{|I|}(\bigodot_{i\in I} x_i)
\right)
\\ &\quad =
\frac{1}{2} \sum_{\substack{I\sqcup J=[n] \\ |I|,|J|\geq 1} }\epsilon(I,J) (-1)^{|k_I|}
\left(
\mathcal{F}_{|I|}(\bigodot_{i\in I} x_i) \circ \mathcal{F}_{|J|}( \bigodot_{j\in J} x_j) 
\right. -\\ &\quad\quad\quad\quad\quad\left. - 
(-1)^{(|k_I|+1)(|k_J|+1)+|k_I||k_J| -|k_I|+|k_J|} \mathcal{F}_{|I|}(\bigodot_{i\in I} x_i) \circ \mathcal{F}_{|J|}( \bigodot_{j\in J} x_j) 
\right)
\\ &\quad =
\sum_{\substack{I\sqcup J=[n] \\ |I|,|J|\geq 1} }\epsilon(I,J) (-1)^{|k_I|}
\mathcal{F}_{|I|}(\bigodot_{i\in I} x_i) \circ \mathcal{F}_{|J|}( \bigodot_{j\in J} x_j) 
\end{align*}

where we use the shorthand $|k_I|=\sum_{i\in I} k_i$ and switched the summation variables $I$ and $J$ for the second equality.
Note that since $d_H=\co{m}{\cdot}_G$, we can absorb the first term of \eref{equ:dglaLinfty} into this expression merely by admitting $I,J=\emptyset$ in the sum and defining $\mathcal{F}_0:= m$.

On the polyvector field side
\begin{align*}
&\frac{1}{2}\sum_{i\neq j}^n \epsilon (i,j,1,..,\hat{i},..,\hat{j},.., n) \mathcal{F}_{n-1}( Q_2(x_i\odot x_j)\odot x_1\odot\cdots\odot \hat{x}_i \odot \cdots \odot\hat{x}_j \odot \cdots x_n)
\\ &\quad=
-\frac{1}{2}
\sum_{i\neq j}^n \epsilon (i,j,1,..,\hat{i},..,\hat{j},.., n) 
\\
&\quad\quad\quad\quad
\mathcal{F}_{n-1}((x_i\bullet x_j+(-1)^{k_i k_j} x_j\bullet x_i)\odot x_1\odot\cdots\odot \hat{x}_i \odot \cdots \odot\hat{x}_j \odot \cdots x_n)
\\ &\quad=
-\sum_{i\neq j}^n \epsilon (i,j,1,..,\hat{i},..,\hat{j},.., n) \mathcal{F}_{n-1}((x_i\bullet x_j)\odot x_1\odot\cdots\odot \hat{x}_i \odot \cdots \odot\hat{x}_j \odot \cdots x_n)
\end{align*}

Hence the conditions \eref{equ:dglaLinfty} for $\mathcal{F}$ to be an $L_\infty$-morphism can be rewritten as
\begin{multline}
\label{equ:ourLinfty}
\sum_{i=1}^n (-1)^{\sum_{r=1}^{i-1}k_r} \mathcal{F}_n( x_1\odot \cdots\odot u \dv x_i \odot \cdots \odot x_n)
-\\-
\sum_{i\neq j}^n \epsilon (i,j,1,..,\hat{i},..,\hat{j},.., n) \mathcal{F}_{n-1}((x_i\bullet x_j)\odot x_1\odot\cdots\odot \hat{x}_i \odot \cdots \odot\hat{x}_j \odot \cdots \odot x_n)
=\\=
\sum_{I\sqcup J=[n]}\epsilon(I,J) (-1)^{|k_I|}
\mathcal{F}_{|I|}(\bigodot_{i\in I} x_i) \circ \mathcal{F}_{|J|}( \bigodot_{j\in J} x_j) 
\end{multline}

\begin{rem}
Note that we can replace all $\odot$'s in the above formula by $\otimes$'s, and the reader should not be worried if this happens soon. In fact, the $\odot$'s are merely a reminder that the functions $\mathcal{F}_n$ on $V^{\otimes n}$ are symmetric, i.e., vanish on the ideal $I$ (intersected with $V^{\otimes n}$).
\end{rem}

\section{Kontsevich and Cyclic Morphism}
\label{sec:graphs}
\subsection{Kontsevich morphism}
In his famous paper \cite{kontsevich} M. Kontsevich constructed an $L_\infty$ quasi-isomorphism 
\[
\morphU : T_{poly}^{\bullet+1}(M) \to D_{poly}^{\bullet+1}(M).
\]
In this subsection we recall his construction for $M=\R^d$, slightly adapted to simplify later proofs. In particular, we will throughout work with the complex $CC^\bullet(M)$ instead of $D^\bullet_{poly}(M)$, which is equivalent due to Lemma \ref{lem:dccisom}. 
 
The morphism can be expressed as a sum of graphs. Denote by $\morphU_m$ the $m$-th Taylor component of $\morphU$. It is given on polyvector fields $\gamma_1,\dots,\gamma_m\in T_{poly}^{\bullet+1}$ by 
\[
\morphU_m(\gamma_1 \otimes \dots \otimes\gamma_m) = [\tilde{\morphU}_m(\gamma_1\otimes\dots \otimes\gamma_m)]
\]
where $\tilde{\morphU}_m(\gamma_1\otimes\dots \otimes\gamma_m)\in \widetilde{CC}^\bullet(M)$ is defined as 
\begin{multline*}
\tilde{\morphU}_m(\gamma_1\otimes\dots \otimes\gamma_m)(a_0\otimes\dots\otimes a_n)
= \sum_{\Gamma\in G(m,n)} w_\Gamma D_\Gamma(\gamma_1\otimes\dots \otimes\gamma_m ; a_0,..,a_n)\omega.
\end{multline*}
The sum is over all Kontsevich graphs with $m$ type I and $n+1$ type II vertices. 

\begin{defi}
\label{def:kontsgraph}
The set $G(m,n)$, $m,n\in \mathbb{N}_0$ of \emph{Kontsevich graphs} consists of directed graphs $\Gamma$ such that
\begin{enumerate}
\item The vertex set of $\Gamma$ is 
\[
V(\Gamma) = \{1,..,m\}\cup \{\bar{0},..,\bar{n}\}
\]
where the vertices $\{1,..,m\}$ will be called the \emph{type I} vertices and the vertices $\{\bar{0},..,\bar{n}\}$ the \emph{type II} vertices.
\item Every edge $e=(v,w)\in E(\Gamma)$ starts at a type I vertex, i.e., $v\in\{1,\dots,m\}$.
\item For each type I vertex $j$, there is an ordering given on
\[
Star(j) = \{(j, w) \mid (j, w)\in E(\Gamma),\; w\in V(\Gamma) \}.
\] 
\item There are no double edges, i.e., edges $(j,w)$ occuring twice in $E(\Gamma)$.
\item There are no tadpoles, i.e., edges of type $(j,j)$.
\end{enumerate}
\end{defi}

The function $\Phi:=D_\Gamma(\gamma_1\otimes\dots \otimes\gamma_m ; a_0,..,a_n)$ is defined as follows.
Let, in standard coordinates on $\R^d$,
\begin{equation}
\label{equ:gammaallind}
\gamma_j = \gamma_j^{i_1\dots i_{k_j}} \p_{i_1}\cdots \p_{i_{k_j}}.
\end{equation}
Here the implicit sum runs over all indices $i_1,..,i_{k_j}=1,..,d$.
Denote by $e_1^j,e_2^j,..$ the edges in $\Gamma$ starting at vertex $j$ in the order as given in the data defining a Kontsevich graph. Let $f_1^v,f_2^v,..$ be the edges ending at vertex $v$ in an arbitrary order. Then
\begin{equation}
\label{equ:dgammadef}
\Phi = \sum_{\varphi:E(\Gamma)\to [d]} \prod_{j=1}^m (\p_{\varphi(f_1^j)} \p_{\varphi(f_2^j)}\cdots \gamma_j^{\varphi(e_1^j)\varphi(e_2^j)\cdots}) \prod_{k=0}^n (\p_{\varphi(f_1^{\bar{k}})} \p_{\varphi(f_2^{\bar{k}})}\cdots a_k)
\end{equation}
where the sum runs over all maps $\varphi$ from the edge set of $\Gamma$ to the set $\{1,..,d\}$.

Let us next define the weight $w_\Gamma$ of $\Gamma\in G(n,m)$. It is an integral of a certain differential form over a compact manifold with corners, the \emph{configuration space} $C_\Gamma$.
\begin{equation}
\label{equ:wgammadef}
w_\Gamma =  \int_{C_\Gamma} \omega_\Gamma
\end{equation}
%Here $k_j=|Star(j)|$ and $k_{[n]}! :=k_1!\cdots k_n!$. \frac{1}{k_{[n]}!} 
\begin{rem}
\label{rem:wedges}
(i) Note that in our conventions the usual factor $\prod_j \frac{1}{|Star(j)|!}$ in the definition of $w_\Gamma$ is missing. To compensate, the sum in \eref{equ:gammaallind} runs over all sets of indices, not just ordered sets.
(ii) Note that due to the summation over orderings of each star, it is unnecessary to require that the tensor $\gamma_j^{i_1\dots i_{k_j}}$ occuring in \eref{equ:gammaallind} is antisymmetric. In fact, one could define $D_\Gamma(\dots)$ also on non-antisymmetric tensors $\gamma_j^{i_1\dots i_{k_j}}$. This fact will be helpful later.
\end{rem}

\begin{defi}
\label{def:configspace}
The {enlarged configuration space} $\tilde{C}_\Gamma$ is the Fulton-MacPherson-like\footnote{We mean the compactification constructed in \cite{kontsevich}, section 5. We will not repeat the construction here.} compactification of the space of embeddings 
\[
(z_1,\dots,z_m,z_{\bar{0}},\dots z_{\bar{n}}): V(\Gamma) \rightarrow D
\]
of the vertex set $V(\Gamma)$ of $\Gamma$ into the closed unit disk $D=\{z \in \cn{}; |z|\leq 1\}$ such that
\begin{enumerate}
\item All type I vertices are mapped to the interior of $D$, i.e. $z_j\in D^\circ$ for $j=1,..,m$.
\item All type II vertices are mapped to the boundary of $D$, i.e. $z_{\bar{j}}\in \partial D$ for $j=0,..,n$.
\item The type II vertices occur in counterclockwise increasing order on the circle, i.e., $0<\arg \frac{z_{\bar{1}}}{z_{\bar{0}}}< \dots < \arg \frac{z_{\bar{n}}}{z_{\bar{0}}}<2\pi$.
\end{enumerate}
The \emph{configuration space} $C_\Gamma$ is the quotient of $\tilde{C}_\Gamma$ under the action of the automorphism group of the unit disk $SU(1,1)$. We put the orientation on $\tilde{C}_\Gamma$ defined by the form 
\[
\Omega=dx_1\wedge dy_1\wedge \cdots\wedge dx_m\wedge dy_m \wedge d\phi_0 \wedge d\phi_n \cdots d\phi_1
\]
where $z_k = x_k+iy_k$ and $z_{\bar{k}} = \exp(2\pi i\phi_k)$. We put on $C_\Gamma$ the induced orientation.\footnote{We mean the orientation on $C_\Gamma$ defined induced by the form $\iota_t\iota_s\iota_h \Omega$ where $h$ is the counterclockwise rotation generator, $s$ the generator of scalings in the upper halfplane model of the hyperbolic disk, and $t$ the generator of right translations of the halfplane.}
\end{defi}

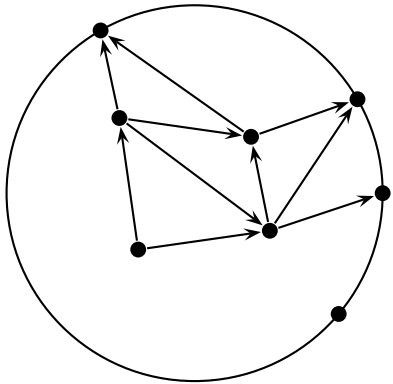
\begin{figure}
\psset{arrowscale=1.5}
\begin{pspicture}(-3,-3)(3,3)
  \pscircle(0,0){2.5}
  \cnode*(-1,1){3pt}{A}
  \cnode*(1,-.5){3pt}{B}
  \cnode*(.75,.75){3pt}{C}
  \cnode*(-.75,-.75){3pt}{D}
  \cnode*(2.5;0){3pt}{x} 
  \cnode*(2.5;30){3pt}{y}
  \cnode*(2.5;120){3pt}{z}
  \cnode*(2.5;-40){3pt}{w}
  
  \psset{arrows=->}
  \ncline{A}{B}
  \ncline{A}{C}
  \ncline{A}{z}
  \ncline{B}{C}
  \ncline{D}{B}
  \ncline{D}{A}
  \ncline{C}{z}
  \ncline{B}{x}
  \ncline{B}{y}
  \ncline{C}{y}
\end{pspicture}
\caption{\label{fig:exgraph} A Kontsevich graph.}
\end{figure}

An example graph embedded in $D$ is shown in Figure \ref{fig:exgraph}.

The differential form $\omega_\Gamma$ that is integrated over configuration space can be expressed as a product of 
one-forms, one for each edge in $\Gamma$. 
\begin{align*}
\omega_\Gamma 
&= \bigwedge_{j=1}^n \bigwedge_{(j, v)\in E(\Gamma)} \alpha(j,v)
\end{align*}

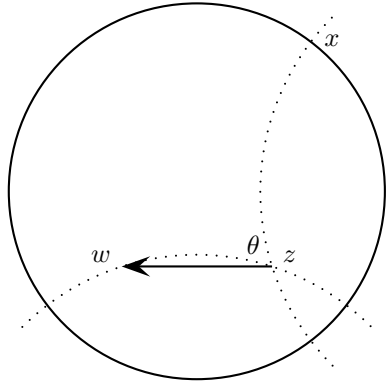
\begin{figure}
\psset{arrowscale=2.5}
\begin{pspicture}(0,-1)(6,6)
  \pscircle(3,3){2.5}
  \pnode(4,2){z}
  \pnode(2,2){w}
  \pnode(4.7,5){x}
  \psline{->}(z)(w) %(4,2)(2,2)
  \psarc[linestyle=dotted](3,-1.125){3.2811}{45}{135}
  \psarc[linestyle=dotted](7.125,3){3.2811}{135}{225}
  \uput{1.5ex}[135](z){$\theta$}
  \uput{1ex}[30](z){$z$}
  \uput{1ex}[150](w){$w$}
  \uput{0}[0](x){$x$}
\end{pspicture}
\caption{\label{fig:phis} Geometric meaning of Kontsevich's angle forms.}
\end{figure}

Here the one-form is defined as $\alpha(j,v)= d\theta(z_j,z_v, z_{\bar{0}})$ where
\begin{align}
\label{equ:thetadef}
d\theta(z,w,x) &= \frac{1}{2 \pi} d\arg \left( \frac{(w-z)(1-\bar{z}x)}{(1-\bar{z}w)(x-z)} \right)
\end{align}
is the differential of the hyperbolic angle between the hyperbolic straight lines $(z,x)$ and $(z,w)$, increasing in the counterclockwise direction (see Figure \ref{fig:phis}). The ordering of the forms within the wedge products is such that forms corresponding to edges with source vertex $j$ stand on the left of those with source vertex $j+1$, and according to the order given on the stars for edges having the same source vertex.

\begin{rem}
As we have written it, all the differential forms above are actually defined on the enlarged configuration space $\tilde{C}_\Gamma$. But one can check that they are $SU(1,1)$-basic and hence descend to the quotient.
\end{rem}

\begin{rem}
\label{rem:thetac}
Note that the form $d\theta(z, w, x)$ satisfies $d\theta(z, w,x) = d\theta(z, w,x') + d\theta(z, x',x)$ for any $x'\in D\setminus\{z\}$. This is very important.
\end{rem}

\subsection{Cyclic Kontsevich morphism}
In this section we define the cyclic variant 
\[
\morphU^{cyc}: \left( T_{poly}(M)[u], u \dv_\omega, \co{}{}_{S} \right) \to \left( (CC(M))^\sigma, d_H, \co{}{}_{G} \right)
\]
of Kontsevich's morphism for $M=\R^d$ and $\omega$ a constant volume form. The Taylor components are $\morphU^{cyc}_n=[\tilde{\morphU}^{cyc}_n]$, where 
\begin{multline*}
\tilde{\morphU}^{cyc}_n(u^{j_1}\gamma_1\otimes \cdots \otimes u^{j_m}\gamma_m)(a_0,\dots, a_n)
=
\sum_{\Gamma\in G_{ex}(m,n)} w_\Gamma(j_1,..,j_m) D_\Gamma(\gamma_1\otimes\dots \otimes\gamma_m ; a_0,..,a_n)\omega.
\end{multline*}
Here, the sum is over all \emph{extended Kontsevich graphs}, which are by definition just Kontsevich graphs with tadpoles.
\begin{defi}
An \emph{extended Kontsevich graph} is a graph satisfying the requirements of Definition \ref{def:kontsgraph}, except possibly the no-tadpole-property (5). We call the set of such graphs with $m$ type I and $n+1$ type II vertices $G_{ex}(m,n)$. For a graph $\Gamma\in G_{ex}(m,n)$, we call the set of vertices with tadpoles $Tp(\Gamma)\subset V(\Gamma)$.
\end{defi}

The poly-differential operator $D_\Gamma$ on the right is defined by exactly the same formula \eref{equ:dgammadef} as in the Kontsevich case. Essentially, this amounts to inserting the divergences of polyvector fields at tadpole vertices, and ignoring the tadpoles otherwise.

Also as before, the weight $w_\Gamma(j_1,..,j_n)$ is computed as an integral over configuration space 
\begin{equation}
\label{equ:wgammacycdef}
w_\Gamma(j_1,..,j_n) = \int_{C_\Gamma} \omega_\Gamma(j_1,..,j_n).
\end{equation}
However, the weight form $\omega_\Gamma(j_1,..,j_n)$ is defined slightly differently, and in particular depends on the $u$-degrees $j_1,\dots,j_m$ of the polyvector fields inserted. Concretely
\begin{align*}
\omega_\Gamma 
&=
\bigwedge_{i=1}^n \varpi^{j_i}_{z_i} \wedge \bigwedge_{(i, v)\in E(\Gamma)} \alpha(i,v).
\end{align*}
with the following definitions:
\begin{itemize}
\item To a tadpole edge, we associate the 1-form $\alpha(i,i) := \eta_{z_i}$.
\item The \emph{non-closed} 1-form $\eta_{z}$ is defined as follows:
\begin{equation}
\label{equ:etadef}
\eta_{z} = \sum_{i=0}^n \theta(z,z_{\overline{i+1}}, z_{\bar{i}}) d\theta(z,z_{\bar{i}}, z_{\bar{0}})
\end{equation}
Here the function $\theta(z,z_{\overline{i+1}}, z_{\bar{i}})$ taking values in $[0,1]$ is defined as in \eref{equ:thetadef}, but with the differentials omitted. It is a well defined smooth function, since both $z_{\overline{i+1}}$ and $z_{\bar{i}})$ lie on the boundary of the disk.

\item The form $\varpi_{z_i}$ is the closed 2-form:
\begin{equation}
\label{equ:varpidef}
\varpi_{z} = -d\eta_{z} =\sum_{i=0}^n d\theta(z,z_{\bar{i}}, z_{\bar{0}})\wedge d\theta(z,z_{\overline{i+1}}, z_{\bar{i}})
\end{equation}
\end{itemize}

Note that the forms $\eta_{z}$ and $\varpi_{z}$ depend on all $z_{\bar{i}}$, though we do not make this dependence explicit to simplify the notation.

We can summarize the above construction of $\morphU^{cyc}$ sloppily by saying that we take Kontsevich's morphism on $T_{poly}(M)$ and extend it to $T_{poly}(M)[u]$ in the following manner:
\begin{enumerate}
\item Replace all $u$'s by $\varpi_z$'s.
\item Allow tadpole graphs and assign the weight forms $\eta_{z}$ to the tadpole edges.
\end{enumerate}

\section{$\morphU^{cyc}$ is cyclically invariant}
\label{sec:cycinv}
The goal of this section is the following Proposition:
\begin{prop}
\label{prop:cycinv}
The pre-$L_\infty$-morphism $\morphU^{cyc}$ constructed in the last section takes values in the cyclically invariant subspace $(CC^\bullet(M))^\sigma\subset CC^\bullet(M)$.
\end{prop}

Having shown this result, we will use the symbol $\morphU^{cyc}$ also for the morphisms $T_{poly}[u] \to (CC^\bullet(M))^\sigma$ and $T_{poly}[u] \to (D_{poly}^\bullet(M))^\sigma$.

\begin{rem}
The special case of the proposition, with all polyvector fields contained in $(T_{poly}(M))_{div}$ was proven in \cite{fsh}.
\end{rem}

\subsection{The cyclic shift}
We defined above the action on the cyclic group of order $n+1$ on cochains generated by
\[
\sigma : \widetilde{CC}^n(M) \to \widetilde{CC}^n(M).
\]
We can also define the action of the cyclic group on the set of graphs $G_{ex}(m,n)$ by cyclically relabeling the type II vertices. Concretely, it is defined such that
\[
\sigma (D_\Gamma(\dots)\omega) = (-1)^n D_{\sigma \Gamma}(\dots)\omega.
\]
with $D_\Gamma(\dots)$ being the polydifferential operator associated to a graph, as in the previous sections. There is also a natural isomorphism between the configuration spaces $\sigma: C_{\sigma \Gamma} \to C_\Gamma$ defined such that
\[
\sigma^* z_{\bar{j}} = z_{\overline{j+1}}
\]
where the $z_{\bar{j}}$ on the left is one of the coordinate functions on $C_\Gamma$, while $z_{\overline{j+1}}$ on the right is a coordinate function on $C_{\sigma \Gamma}$.\footnote{This notation is not optimal. The different $z_{\bar{j}}$'s for different configuration spaces should be distinguishable.}

Using these definitions, we can compute
\[
\sigma (\sum_\Gamma w_\Gamma D_\Gamma\omega)
=
\sum_\Gamma (-1)^n w_\Gamma D_{\sigma\Gamma}\omega
=\sum_\Gamma (-1)^n w_{\sigma^{-1}\Gamma} D_{\Gamma}\omega
\]
Here and in the remainder of this section we suppress the arguments $j_1,..,j_m$ of $w_\Gamma$ and $\gamma_1,..,\gamma_m$ of $D_\Gamma$ since they do not play any role. Our goal is to show that 
\[
\sum_{\Gamma\in G_{ex}(m,n)} (-1)^n w_{\sigma^{-1}\Gamma} [D_{\Gamma}\omega]
=
\sum_{\Gamma\in G_{ex}(m,n)} w_{\Gamma} [D_{\Gamma}\omega]
\]

\subsection{First Simplification}
\begin{defi}
A \emph{marked graph} is an extended Kontsevich graph $\Gamma$ together with a distinguished subset of its edges (``marked edges''). 
We denote the sets of such graphs by $G_{ex}'(m,n)$. To a marked graph $\Gamma\in G_{ex}'(m,n)$ we associate a weight $w_\Gamma$ by formulas as above, but with the weight forms $\alpha(j,v) = -\alpha(j,\bar{1})$ associated to its marked edges.
\end{defi}

\begin{lemma}
\label{lem:markedgr}
Let $\Gamma \in G_{ex}(m,n)$. Then
\[
(-1)^n w_{\sigma^{-1}\Gamma} = \sum_{\Gamma'}w_{\Gamma'} 
\]
where the sum on the right runs over all marked graphs $\Gamma'\in G_{ex}'(m,n)$ whose underlying extended Kontsevich graph is $\Gamma$, i.e., over all ``markings'' of $\Gamma$.
\end{lemma}

This result will be an easy consequence of the following lemma:

\begin{lemma}
\label{lem:angleforms}
Let $\Gamma\in G_{ex}(m,n)$ and $\sigma: C_{\Gamma} \to C_{\sigma^{-1}\Gamma}$ as above.
Denote the coordinate functions on $C_\Gamma$ by $z_{\bar j}$ and on $C_{\sigma^{-1} \Gamma}$ by $z_{\bar j}'$, and similarily distinguish $\varpi_z$ on $C_\Gamma$ from $\varpi_z'$ on $C_{\sigma^{-1}\Gamma}$ etc. Then
\begin{enumerate}
\item $\sigma^* \varpi_z'= \varpi_z$.
\item $\sigma^* \alpha(i,v)' = \alpha(i,\sigma(v))-\alpha(i,\bar{1})$, where $\sigma(v)=v$ for $v$ of type I and $\sigma(\bar{j})=\overline{\sigma(j)}$.
\end{enumerate}
\end{lemma}

In words, the Lemma says that the forms $\varpi_z$ are invariant under cyclic interchange of the type II vertices, and that the forms associated to edges $(i,v)$ change by the same form that would be associated to an edge $(i,\bar{1})$.

\begin{proof}
(1) For brevity, define $u_j= \theta(z,z_{\overline{j}}, z_{\bar{0}})$ and $\sigma_j=\theta(z,z_{\overline{j+1}}, z_{\bar{j}})$, and similarly $u_j'$ and $\sigma_j'$. Using Remark \ref{rem:thetac}, it is easily seen that 
\begin{align*}
u_j &= \sum_{i=0}^{j-1} \sigma_i &\quad \quad
\sum_{i=0}^{n} \sigma_i &= 1 \\
\varpi_z &= \sum_{i=0}^n du_i \wedge du_{i+1} &\quad\quad
\sigma^* \sigma_j' &= \sigma_{j+1}  \\
\sigma^* u_j' &= u_{j+1} - \sigma_{0}
\end{align*}

Compute
\begin{align*}
\sigma^* \varpi_z' &= \sigma^* \sum_{i=0}^n du_i' \wedge du_{i+1}'
= \sum_{i=0}^n d(u_{i+1}-\sigma_0) \wedge d(u_{i+2}-\sigma_0) \\
&= \sum_{i=0}^n du_{i+1} \wedge du_{i+2} 
-\left( \sum_{i=0}^n du_{i+1}\right) \wedge d\sigma_0
-d\sigma_0\wedge \left(\sum_{i=0}^n du_{i+2}\right) \\
&= \varpi_z + 0.
\end{align*}

(2) We distinguish the cases $i=v$ (tadpole) and $i\neq v$ (no tadpole). In the latter case, the statement immediately follows from Remark \ref{rem:thetac}. Concretely,
\[
\sigma^*d\theta(z_i',z_v',z_{\bar{0}}') =  d\theta(z_i,z_{\sigma(v)},z_{\bar{1}})
= d\theta(z_i,z_{\sigma(v)},z_{\bar{0}}) + d\theta(z_i,z_{\bar{0}},z_{\bar{1}})
= \alpha(i,\sigma(v)) - \alpha(i, \bar{1}).
\]
 
For the former (tadpole) case, note that
\[
\eta_z = \sum_{i=0}^n \sigma_i du_{i}
\]
and hence, similarly to the calculation above
\begin{align*}
\sigma^* \eta_z' &= \sum_{i=0}^n \sigma_{i+1} d(u_{i+1}-\sigma_0) 
= \eta_z - \left( \sum_{i=0}^n \sigma_{i+1} \right) d\sigma_0 \\
&= \eta_z - d\sigma_0.
\end{align*}
\end{proof}

\begin{proof}[Proof of Lemma \ref{lem:markedgr}]
Note that $\sigma$ is an isomorphism of the configuration spaces that changes the orientation by a factor $(-1)^n$. Hence
\begin{align*}
(-1)^n w_{\sigma^{-1}\Gamma} &= (-1)^n \int_{C_{\sigma^{-1}\Gamma}}\omega_\Gamma \\
&= \int_{C_\Gamma}\sigma^*\omega_{\sigma^{-1}\Gamma} 
&= \int_{C_\Gamma} \bigwedge_{i=1}^n \varpi^{j_i}_{z_i} \wedge \bigwedge_{(i, v)\in E(\Gamma)} (\alpha(i,v)-\alpha(i,\bar{1}))
&= \sum_{\Gamma'}w_{\Gamma'}
\end{align*}
Here we used Lemma \ref{lem:angleforms} in the third equality.
\end{proof}

\subsection{Proof of Cyclic Invariance}
The following lemma is essentially ``integration by parts''.
\begin{lemma}
\label{lem:dottededges}
Let $\Gamma\in G_{ex}(m,n)$ with an edge $(j,\bar{0})$.
The cochain $[D_\Gamma\omega] \in CC^\bullet(M)$ equals a sum of cochains $-\sum_{\Gamma'} [D_{\Gamma'}\omega]$, where the sum runs over all graphs $\Gamma'$ obtained by replacing the edge $(j,\bar{0})$ with the edge $(j,v)$ for $v$ an arbitrary vertex in $\Gamma$ except $\bar{0}$.
\end{lemma}
\qed

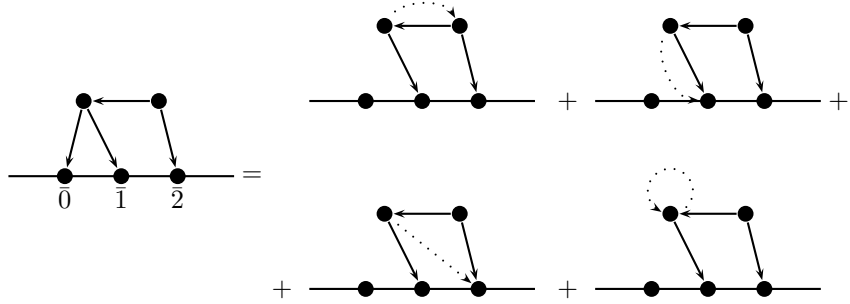
\begin{figure}
\psset{unit=1cm}
\psset{arrowscale=1}
\psset{coilwidth=0.2, coilarm=0.25, coilaspect=0,coilheight=2}
%\begin{tabular}{r c l}
%\multirow{2}{*}{
%\centering
%\begin{pspicture}(0,-.5)(3,1.5)
\begin{pspicture}(0,-2)(11,3)  
  \cnode*(1,1){3pt}{A}
  \cnode*(2,1){3pt}{B}
  \cnode*(.75,0){3pt}{x0}
  \cnode*(1.5,0){3pt}{x1}
  \cnode*(2.25,0){3pt}{x2}   
  \psline{-}(0,0)(3,0)
  \psset{arrows=->}
  \ncline{B}{A}
  \ncline{A}{x0}
  \ncline{A}{x1}
  \ncline{B}{x2}
  \uput{1ex}[-90](x0){$\bar{0}$}
  \uput{1ex}[-90](x1){$\bar{1}$}
  \uput{1ex}[-90](x2){$\bar{2}$}
  \uput{0}[0](3.1,0){$=$}
  \uput{0}[0](7.3,1){$+$}
  \uput{0}[0](7.3,-1.5){$+$}
  \uput{0}[0](3.5,-1.5){$+$}
  \uput{0}[0](10.9,1){$+$}
  
  \cnode*(5,2){3pt}{A1}
  \cnode*(6,2){3pt}{Bb}
  \cnode*(4.75,1){3pt}{x01}
  \cnode*(5.5,1){3pt}{x11}
  \cnode*(6.25,1){3pt}{x21} 
  
  \psline{-}(4,1)(7,1)
  \ncline{Bb}{A1}
  \ncarc[arcangle=60, linestyle=dotted]{A1}{Bb}
  \ncline{A1}{x11}
  \ncline{Bb}{x21}

  \cnode*(8.8,2){3pt}{A2}
  \cnode*(9.8,2){3pt}{B2}
  \cnode*(8.55,1){3pt}{x02}
  \cnode*(9.3,1){3pt}{x12}
  \cnode*(10.05,1){3pt}{x22} 
  
  \psline{-}(7.8,1)(10.8,1)
  \ncline{B2}{A2}
  \ncline{A2}{x12}
  \ncline{B2}{x22}
  \ncarc[arcangle=-60, linestyle=dotted]{A2}{x12}

  \cnode*(5,-.5){3pt}{A3}
  \cnode*(6,-.5){3pt}{B3}
  \cnode*(4.75,-1.5){3pt}{x03}
  \cnode*(5.5,-1.5){3pt}{x13}
  \cnode*(6.25,-1.5){3pt}{x23} 
  
  \psline{-}(4,-1.5)(7,-1.5)
  \ncline{B3}{A3}
  \ncline{B3}{x23}
  \ncline{A3}{x13}
  \ncline[linestyle=dotted]{A3}{x23}

  \cnode*(8.8,-.5){3pt}{A4}
  \cnode*(9.8,-.5){3pt}{B4}
  \cnode*(8.55,-1.5){3pt}{x04}
  \cnode*(9.3,-1.5){3pt}{x14}
  \cnode*(10.05,-1.5){3pt}{x24} 
  
  \psline{-}(7.8,-1.5)(10.8,-1.5)
  \ncline{B4}{A4}
  \ncline{B4}{x24}
  \ncline{A4}{x14}
  \nccircle[linestyle=dotted]{->}{A4}{.3}
\end{pspicture}
%\end{tabular}
\caption{\label{fig:dottededges}Graphical illustration of the meaning of Lemma \ref{lem:dottededges}.}
\end{figure}

\begin{proof}[Proof of Proposition \ref{prop:cycinv}]
Given Lemma \ref{lem:markedgr}, we need to show that 
\begin{equation}
\label{equ:cycinvtbs}
\sum_{\Gamma\in G_{ex}(m,n)} w_{\Gamma}(\dots) [D_{\sigma\Gamma}(\dots)\omega]
=
\sum_{\Gamma\in G'_{ex}(m,n)} w_{\Gamma}(\dots) [D_{\sigma\Gamma}(\dots)\omega].
\end{equation}
Apply Lemma \ref{lem:dottededges} to each marked edge hitting $\bar{0}$ of the marked graphs on the right hand side. Graphically we remove each marked edge $(j,\bar{0})$ and reconnect it to all other vertices except $\bar{0}$, see Figure \ref{fig:dottededges}. Let us call the reconnected edges ``dotted edges''. Assign to these dotted edges the weight form $\alpha(j,\bar{1})$, i.e., the negative of the weight form of the marked edge we removed. The right hand side of \eref{equ:cycinvtbs} can hence be rewritten as
\[
\sum_{\Gamma} w_{\Gamma}(\dots) [D_{\sigma\Gamma}(\dots)\omega]
\]
where the sum runs over all marked graphs without marked edges connecting to $\bar{0}$, but instead with arbitrarily many additional dotted edges. This sum can be simplified as follows.

\textbf{First Cancellation:} For a graph $\Gamma$ which contains double edges, $D_{\Gamma}$ vanishes. This is because polyvector fields are antisymmetric.

\textbf{Second Cancellation:} The contribution of any graph $\Gamma$ containing a dotted edge $(j,v)$ cancels with the contribution of the graph $\Gamma'$, in which the dotted edge $(j,v)$ is replaced by a marked edge $(j,v)$. This is because the weight forms of these edges sum up to 0.

Hence all graphs with marked and/or dotted edges cancel each other. What remains are the graphs containing only solid unmarked edges, i.e., exactly those on the left of \eref{equ:cycinvtbs}.
\end{proof}

\section{The proof of Theorem \ref{thm:main}}
\label{sec:theproof}

In this section we will prove the main Theorem \ref{thm:main}. There are two things left to be shown:
\begin{enumerate}
\item The morphism $\morphU^{cyc}$ is a quasi-isomorphism. This will follow from Shoikhet's work and the fact that on divergence free vector fields our morphism agrees with Shoikhet's.

\item It is an $L_\infty$-morphism of dglas. The proof is a typical ``Kontsevich-Stokes-proof''. The only unusual thing is that we have a non-vanishing ``bulk'' term due to the non-closedness of $\eta_z$, which will provide the differential $u \dv$ on the polyvector field side. This trick has probably been invented by Cattaneo and Felder \cite{cafe}.
\end{enumerate}

\subsection{It is a Quasi-Isomorphism}
This is a direct consequence of the following Proposition and Theorem 2.3 of \cite{shoikhet}:
\begin{prop}
On the sub-Lie algebra $(T_{poly}(M))_{div}[u]$ of divergence free poly vector fields, the morphism $\morphU^{cyc}$ agrees with Shoikhet's morphism from \cite{shoikhet}, up to signs due to different conventions.
\end{prop}
\begin{proof}
Shoikhet's morphism contains a sum of graphs with dashed pairs. Using the notation of the proof of Lemma \ref{lem:angleforms}, a dashed pair at $i$-th position introduces a factor $du_i\wedge du_{i+1}$ into the weight form. Summing over all $i$ we get just the form $\varpi_z$. 
\end{proof}

% \begin{rem}
% A more thorough analysis using the above methods shows that the full Shoikhet morphism can be constructed along the same lines as ours, but with the single replacement $\eta_z\to +-??? d\theta(z, z_{\bar{1}},z_{\bar{0}})$ ???. 
% \end{rem}

\subsection{Quadratic relations between weights}
\label{sec:qweights}
The proof of $\morphU^{cyc}$ being an $L_\infty$-morphism will closely follow Kontsevich's proof of his formality theorem. In fact, we will sometimes allow ourselves to be sketchy and point out only the differences. Kontsevich's idea was to derive quadratic relations between weights of graphs using Stokes' Theorem :
\begin{equation}
\label{equ:stokes}
\int_{C_\Gamma} d\omega_\Gamma(j_1,\dots ,j_m) = \int_{\p C_\Gamma} \omega_\Gamma(j_1,\dots ,j_m).
\end{equation}

In the Kontsevich case, the left hand side is always zero due to closedness of the angle forms. In our case however, the left hand side can be nonzero, if the graph $\Gamma$ contains a tadpole. More precisely, the left hand side equals 
\begin{align*}
\int_{C_\Gamma} d\omega_\Gamma(j_1,..,j_n) 
& =
-\sum_{i\in Tp(\Gamma)} (-1)^{\sum_{r=0}^{i-1}k_r +s(i,i) }\int_{C_\Gamma} \omega_{\Gamma-(i,i)}(j_1,..,j_i+1,..,j_m)
\\ &=
-\sum_{i\in Tp(\Gamma)}(-1)^{\sum_{r=0}^{i-1}k_r +s(i,i) } w_{\Gamma-(i,i)}(j_1,..,j_i+1,..,j_m)
\end{align*}
where the sum is over all vertices $i$ that have a tadpole, and the graph $\Gamma-(i,i)$ is the graph obtained by removing the tapdpole at $i$. The number $s(i,j)$ is the position of the edge $(i,j)$ in the ordering on $Star(i)$, minus 1. For example, if the edge $(i,i)$ is first in the ordering, then $s(i,i)=0$. The sign in front of the sum over $i$ is due to the fact that $d\eta_z=-\varpi_z$, not ``+''. The numbers $k_r$ equal $|Star(r)|$.
Note also that configuration spaces do not depend on the edges and hence $C_\Gamma=C_{\Gamma-(i,i)}$.

Now consider the right hand side of \eref{equ:stokes}. As in the Kontsevich case, the codimension one boundary strata of $C_\Gamma$ are indexed by ``good'' subsets of the vertex set of $\Gamma$, which collapse to as point. A good subset is one such that (i) all its type II vertices are adjacent and (ii) it does not contain $\bar{0}$. In dealing with these subsets, we will fix the following notations:
\begin{enumerate}
\item The set of type I vertices in the collapsing subset is $J\subset [m]$ and $m_2=|J|$. Its complement is denoted $I:=[m]\setminus J$ and $m_1:=m-m_2=|I|$.
\item The integers $l,n_2$ are such that the set of collapsing type II vertices is $\overline{l},..,\overline{l+n_2-1}$. Let further $n_1:=n-n_2+1$.
\item $\Gamma'\subset \Gamma$ is the full subgraph formed by the collapsing vertices. Furthermore $\Gamma/\Gamma$ is the graph obtained from $\Gamma$ by collapsing $\Gamma'$ to one vertex.
\item $j_I$ is shorthand for $j_{i_1},j_{i_2},\dots,j_{i_{m_1}}$ where $I=\{j_{i_1},\dots,j_{i_{m_1}}\}$ and $j_{i_1}<\cdots<j_{i_{m_1}}$. Similarly for $j_J$.
\end{enumerate}

Using these notations, the corresponding boundary stratum has the form $C_{\Gamma/\Gamma'}\times C_{\Gamma'}$.\footnote{For the definition of the configuration space of a graph without type II vertices, we refer to \cite{kontsevich}. This is the only place where we need it.} We denote the projections onto the left and right factors by $\pi_1$ and $\pi_2$. Consider two cases separately:

(A) $n_2=0$. This boundary stratum corresponds to $m_2\geq 2$ points approaching each other, away from the boundary of the disk. Note that for $z$ one of the points collapsing, the forms $\varpi_z$ and $\eta_z$ are basic wrt. the projection $\pi_1$. Hence they do not spoil Kontsevich's argument (Lemma 6.6 in \cite{kontsevich}) that the contribution of such a stratum vanishes unless $m_2=2$. In this case the two vertices have to be connected by exactly one edge by dimensionality reasons, and the integral over $C_{\Gamma'}$ is $1$. Introducing new notation, the integral over this stratum yields a contribution
\[
- \epsilon(i,j,1,..,\hat{i},..,\hat{j},..,m) (-1)^{s(i,j)} w_{\tilde{\Gamma}}(j_i+j_j,j_1,\dots,\hat{j}_i,\dots,\hat{j}_j,\dots,j_m)
\]
where $i,j$ are the two (simply) connected vertices collapsing and $\tilde{\Gamma}$ is the graph obtained by (i) renumbering the vertices such that vertices $i$ and $j$ become $1$ and $2$ and (ii) contracting the edge $(1,2)$. The sign $\epsilon(i,j,1,..,\hat{i},..,\hat{j},..,m)$ coming from the renumbering is defined similar to that in eqn. \eref{equ:dglaLinfty}. The ordering on $Star(1)$ of $\tilde{\Gamma}$ is such that the edges coming from vertex $i$ of $\Gamma$ stand before those coming from vertex $j$.
The sign in front comes from the orientations of the spaces involved, see \cite{arnal}, Lemma I.2.1. 
%In particular case of, if both approaching vertices have a tadpole, $\Gamma / (i,j)$ has two and hence the weight vanishes.

(B) $n_2>0$. This boundary stratum corresponds to $m_2\geq 0$ points in the interior and $n_2\geq 1$ points on the boundary approaching each other. Let us denote the inclusion  $C_{\Gamma/\Gamma'}\times C_{\Gamma'}\hookrightarrow C_\Gamma$ by $\iota$. It was shown in \cite{arnal}, Lemma I.2.2, that $\iota$ changes orientations by a factor (translated into our nomenclature, sign and orientation conventions) $|\iota|=(-1)^{(l+1)(n_2+1)+n_2+n}$. The differential form $\omega_\Gamma(j_1,..,j_m)$ is equal to the form $\epsilon(I,J)\omega_{\Gamma_{I,J}}(j_I,j_J)$, where $\Gamma_{I,J}$ is the graph obtained by renumbering the type I vertices such that those in $I$ stand to the left of those in $J$. The sign $\epsilon(I,J)$ is again defined as in section \ref{sec:linfty}.

We next claim that $\iota^*\omega_{\tilde{\Gamma}}(j_I,j_J) = \pi_1^*\omega_{\Gamma/\Gamma'}(j_I) \wedge \pi_2^* \omega_{\Gamma'}(j_J)$. 
The novelty here in comparison to the Kontsevich case is the possible occurrence of forms $\varpi_z$ and $\eta_z$. Assume first that the point $z$ is not collapsing to the boundary and recall that 
\[
\varpi_z = \sum_{0\leq i<j\leq n} d\sigma_i\wedge d\sigma_j
\]
in the notation of the proof of Lemma \ref{lem:angleforms}. When the vertices $\overline{l},..,\overline{l+n_2-1}$ collapse, the angles $\sigma_l,..,\sigma_{l+n_2-2}$ go to zero and drop out of the sum. The remaining terms are exactly the $\varpi_z$ of $\Gamma/\Gamma'$. A similar arguments holds for
\[
\eta_z = \sum_{0\leq i<j\leq n} \sigma_j d\sigma_i.
\]

Next suppose that $z$ is one of the vertices in $J$ approaching the boundary. Since all $\sigma_j$ are $SU(1,1)$-invariant, we can as well suppose that the vertices of $\Gamma'$ do not collapse, but instead the complement, i.e. all vertices in $I$ and the type II vertices $\overline{l+n_2},..,\bar{n},\bar{0},..,\overline{l-1}$ collapse to a point on the boundary. Then the same argument as before shows that $\varpi_z$ and $\eta_z$ become the $\varpi_z$ and $\eta_z$ of $\Gamma'$.

Hence we compute
\begin{multline*}
\int_{\p_{\Gamma'} C_\Gamma}\omega_\Gamma(j_1,..,j_m)
=
|\iota| \int_{C_{\Gamma/\Gamma'}\times C_{\Gamma'}} \iota^* \omega_\Gamma(j_1,..,j_m)
=\\=
(-1)^{(l+1)(n_2+1)+n_2+n}\epsilon(I,J) \left( \int_{C_{\Gamma/\Gamma'}} \omega_{\Gamma/\Gamma'}(j_I) \right)\left(\int_{C_{\Gamma'}} \omega_{\Gamma'}(j_J)\right).
\end{multline*}

Let us summarize the result of this subsection.
\begin{prop}
\label{prop:weightrel}
Let $\Gamma$ be a Kontsevich graph with $m$ type I vertices and $j_1,..,j_m=0,1,..$ be integers. Then, with the notations from above
\begin{multline*}
-\sum_{i\in Tp(\Gamma)}(-1)^{\sum_{r=0}^{i-1}k_r +s(i,i) } w_{\Gamma-(i,i)}(j_1,..,j_i+1,..,j_m)
= \\=
- \sum_{ \substack{(i,j)\in E(\Gamma) \\ (j,i)\notin E(\Gamma)} } 
(-1)^{s(i,j)} 
\epsilon(i,j,1,..,\hat{i},..,\hat{j},..,m)
w_{\tilde{\Gamma}}(j_i+j_j,j_1,..,\hat{j}_i,..,\hat{j}_j,..,j_m)
+ \\ +
\sum_{\Gamma'\subset \Gamma} (-1)^{(l+1)(n_2+1)+n_2+n} \epsilon(I,J) w_{\Gamma/\Gamma'}(j_I) w_{\Gamma'} (j_J).
\end{multline*}
\end{prop}

\subsection{It is an $L_\infty$-morphism}
We want to show that \eref{equ:ourLinfty} holds for $\mathcal{F}=\morphU^{cyc}$.
Each of the three terms occuring has a representation in terms of graphs. The first term on the left can be identified with the sum
\begin{multline*}
\sum_\Gamma \sum_{i=1}^m (-1)^{\sum_{r=1}^{i-1}k_r} w_{\Gamma}(j_1,..,j_i+1,..,j_m) D_{\Gamma} (\gamma_1 \otimes \cdots \otimes\dv \gamma_i \otimes \cdots \otimes \gamma_m)
=\\=
\sum_\Gamma \sum_{i=1}^m (-1)^{\sum_{r=1}^{i-1}k_r} w_{\Gamma}(j_1,..,j_i+1,..,j_m) \sum_{s(i,i)} (-1)^{s(i,i)} D_{\Gamma+(i,i)} (\gamma_1 \otimes \cdots \otimes \gamma_m)
=\\=
\sum_\Gamma \sum_{i\in Tp(\Gamma)} (-1)^{\sum_{r=1}^{i-1}k_r+s(i,i)} w_{\Gamma-(i,i)}(j_1,..,j_i+1,..,j_m) D_\Gamma(\gamma_1 \otimes \cdots \otimes \gamma_m).
\end{multline*}
Here we substituted $x_i=u^{j_i}\gamma_i$ into \eref{equ:ourLinfty}. In the second line an edge $(i,i)$ is added to $\Gamma$ and put at the $s(i,i)+1$st position in the ordering on $Star(i)$. The sum on the right is over all such possible positions, $s(i,i)=0,..,|Star(i)|-1$. For the last line, we changed the summation variable $\Gamma$.

The second term on the left of \eref{equ:ourLinfty} can be identified with
\begin{multline*}
\sum_\Gamma 
\sum_{i\neq j}^n \epsilon (i,j,..,\hat{i},..,\hat{j},.., m)
w_\Gamma(j_i+j_j,j_1,..,\hat{j}_i,..,\hat{j}_j,..,j_m) \\
D_\Gamma((\gamma_i\bullet \gamma_j)\otimes \gamma_1\otimes\cdots\otimes \hat{\gamma}_i \otimes \cdots \otimes\hat{\gamma}_j \otimes \cdots \gamma_n)
=\\=
\sum_\Gamma
%\sum_{\substack{(i,j)\in E(\Gamma) \\ (j,i)\notin E(\Gamma)}}
\sum_{i\neq j}^n \epsilon (i,j,..,\hat{i},..,\hat{j},.., m) (-1)^{k_i}
w_{\Gamma}(j_i+j_j,j_1,..,..,\hat{j}_j,..,j_m)
\sum_{\Gamma'} (-1)^{s(i,j)} D_{\Gamma'}(\gamma_1\otimes \cdots \otimes \gamma_m)
=\\=
\sum_\Gamma \sum_{(i,j)\in E(\Gamma) } 
(-1)^{s(i,j)} \epsilon (i,j,..,\hat{i},..,\hat{j},.., m)
w_{\tilde{\Gamma}}(j_i+j_j,j_1,..,\hat{j}_i,..,\hat{j}_j,\dots,j_m)
D_\Gamma(\gamma_1\otimes \cdots \otimes \gamma_m).
=\\=
\sum_\Gamma \sum_{\substack{(i,j)\in E(\Gamma) \\ (j,i)\notin E(\Gamma)}} 
(-1)^{s(i,j)} \epsilon (i,j,..,\hat{i},..,\hat{j},.., m)
w_{\tilde{\Gamma}}(j_i+j_j,j_1,..,\hat{j}_i,..,\hat{j}_j,\dots,j_m)
D_\Gamma(\gamma_1\otimes \cdots \otimes \gamma_m).
\end{multline*}
The sum over $\Gamma'$ is over graphs obtained from $\Gamma$ by the following procedure:
\begin{enumerate}
\item Insert an additional type I vertex, and renumber the vertices such that (i) the new vertex is vertex $j$ and (ii) the vertex of $\Gamma$ in which $\gamma_k$ was inserted is the new vertex $1$ in $\Gamma'$. Since the numbering of vertices is irrelevant for the definition of $D_{\Gamma'}$, there does not occur an additional sign.
\item Reconnect zero or more edges ending at vertex $i$ to vertex $j$, i.e., apply Leibniz' rule.  There does not occur a sign either. 
\item Reconnect the $k_j$ last (in the ordering on $Star(i)$) edges starting at $i$ such that they start at $j$, maintaining their ordering. 
\item Add an extra edge $(i,j)$. Make it the $s(i,j)+1$st in the ordering on $Star(i)$, where $s(i,j)=0,1,..,|Star(i)|-1$.
\end{enumerate}
To see the first equality, note that by Remark \ref{rem:wedges} (ii) we are allowed to replace the antisymmetric tensor $\gamma_i\bullet \gamma_j$ by its non-antisymmetrized version $\iota_{\p_r} \gamma_i \otimes \p_r\cdot\gamma_j$. The resulting term
\[
D_\Gamma((\iota_{\p_r} \gamma_i \otimes \p_r\cdot\gamma_j)\otimes \gamma_1\otimes\cdots\otimes \hat{\gamma}_i \otimes \cdots \otimes\hat{\gamma}_j \otimes \cdots \gamma_n)
\]
is equal to the sum
\[
\sum_{\Gamma'} (-1)^{s(i,j)} D_{\Gamma'}(\gamma_1\otimes \cdots \otimes \gamma_m).
\]

For the second equality, we changed the summation variables. Here $\tilde{\Gamma}$ is the graph obtained from $\Gamma$ by renumbering the vertices such that vertices $i$ and $j$ become vertices 1 and 2, maintaining the order of the other vertices. Note that we dropped here graphs containing a double edge $(i,j)$ since they do not contribute due to antisymmetry of the $\gamma$'s. Note that in this sum, there may be graphs containing an edge $(i,j)$ as well as an edge $(j,i)$. More precisely, these graphs come from applying the Leibniz rule to a tadpole edge. For the third equality, we used that all those graphs cancel. Concretely, if there is such a pair of edges, the $(i,j)$-term in the sum over edges cancels with the $(i',j')=(j,i)$-term. The relative sign of the two terms is a $(-1)^{k_ik_j}$ from the $\epsilon(\dots)$, times a $(-1)^{(k_j-1)k_i+k_i-1}$ from permuting the weight forms of the edges of vertices $i$ and $j$ appropriately. Hence the total relative sign is $-1$. This computation is reminiscent of the calculation showing that $\dv$ intertwines with the Schouten bracket.

The term on the right of \eref{equ:ourLinfty} is 
\begin{multline*}
\sum_{I\sqcup J=[n]}\epsilon(I,J) (-1)^{|k_I|} \sum_{\Gamma_1,\Gamma_2} w_{\Gamma_1}(j_I) w_{\Gamma_2}(j_J) D_{\Gamma_1}(\gamma_I) \circ D_{\Gamma_2}(\gamma_J)
=\\=
\sum_{I\sqcup J=[n]}\epsilon(I,J) (-1)^{|k_I|} \sum_{\Gamma_1,\Gamma_2} w_{\Gamma_1}(j_I) w_{\Gamma_2}(j_J) \sum_{l=1}^{n_1} (-1)^{(l+1)(n_2+1)} D_{\Gamma_1}(\gamma_I) \circ_l D_{\Gamma_2}(\gamma_J)
=\\=
\sum_{I\sqcup J=[n]}\epsilon(I,J) (-1)^{|k_I|}
\sum_{\Gamma_1,\Gamma_2} w_{\Gamma_1} w_{\Gamma_2} \sum_{\Gamma_{I,J}}(-1)^{(l+1)(n_2+1)} D_{\Gamma_{I,J}}(\gamma_1\otimes\cdots\otimes\gamma_n)
=\\=
\sum_\Gamma \sum_{\Gamma'\subset \Gamma}\epsilon(I,J) (-1)^{|k_I|+(l+1)(n_2+1)}
w_{\Gamma/\Gamma'}(j_I) w_{\Gamma'} (j_J)
D_\Gamma(\gamma_1\otimes \cdots \otimes \gamma_m).
\end{multline*}

Here the sum over $\Gamma_{I,J}$ in the second line is over all graphs that can be obtained by inserting the graph $\Gamma_2$ into one of the type II vertices $\bar{1},\dots,\overline{n_1}$ of $\Gamma_1$ and reconnecting the edges ending at that type two vertex in any possible way to vertices of $\Gamma_2$. The number $l$ in the second line is such the type II vertex $\Gamma_2$ is substituted for is $\bar{l}$.
In the last line, we again switched summation variables. The numbers $l$ and $n_2$ here are as in Section \ref{sec:qweights}.
To proceed further, note that 
\[
(-1)^{|k_I|+(l+1)(n_2+1)}=(-1)^{n-n_2+1+(l+1)(n_2+1)} = -(-1)^{(l+1)(n_2+1)+n-n_2}
\]
Comparing this to the quadratic weight relations of Proposition \ref{prop:weightrel}, one sees that the three summands sum up to 0. 
Hence Theorem \ref{thm:main} is proven. \qed

\section{Globalization--by Damien Calaque}\label{sec:cycproof}

\subsection{Properties of $\mathcal U^{cyc}$}

Before proving the cyclic formality conjecture for an arbitrary smooth manifold $M$ in the next Subsection, 
we need to prove a certain number of properties of $\mathcal U^{cyc}$ suitable for the globalization. 
\begin{prop}\label{prop:properties}
The $L_\infty$-quasi-isomorphism $\mathcal{U}^{cyc}$ has the following properties: 
\begin{itemize}
\item it can be defined for $\R^d_{\rm formal}$ as well; 
\item $\mathcal{U}^{cyc}_1(\gamma)=\gamma$ for any divergence free 
vector field $\gamma\in(T_{\rm poly}^1(\R^d))_{\rm div}$;
\item $\mathcal{U}_n^{cyc}(\gamma_1,\dots,\gamma_n)=0$ for any $n\geq2$ and any vector fields 
$\gamma_1,\dots,\gamma_n\in T_{\rm poly}^1(\R^d)$; 
\item $\mathcal{U}_n^{cyc}(\gamma,\alpha_2,\dots,\alpha_n)=0$ for any $n\geq2$, 
any divergence free linear vector field 
$\gamma\in\mathfrak{sl}_d(\R)\subset T_{\rm poly}^1(\R^d)$ 
and any elements $\alpha_2,\dots,\alpha_n\in T_{\rm poly}(\R^d)[u]$. 
\end{itemize}
\end{prop}
\begin{proof}
The first property is immediate from the definition of $\mathcal U^{cyc}$. 

The second property follows from the fact that $\gamma$ is divergence free. Therefore 
the graph consisting of a single tadpole does not contribute and the only remaining graph is the one with 
a single edge, going from the type only I vertex to the type II vertex $\overline{1}$. 

For degree reason, the third property is non-trivial only for $n=2$. Then observe that the weight of a graph with 
only one vertex $\overline{0}$ of the second type and at least one tadpole is zero (since in this case the form 
$\eta$ vanishes). Therefore $\mathcal U^{cyc}_2(\gamma_1,\gamma_2)=\mathcal{U}_2(\gamma_1,\gamma_2)=0$ ($\mathcal{U}$ 
being Kontsevich's $L_\infty$-morphism). 

The last property follows form the fact that $\gamma$ is divergence free (therefore there is no tadpole attached to 
its corresponding vertex) and from Kontsevich's vanishing lemma for vertices with exactly one incoming and one 
outgoing edges. 
\end{proof}
%First of all observe that all those properties are satisfied by $\mathcal U^{sh}$. Namely, the proof of 
%the third, resp.~fourth, property is the same as the one of property (P4), resp.~(P5), in \S7.3.1.1, 
%resp.~\S7.3.3.1, of \cite{kontsevich}. And the first two properties are obvious from 
%the definition of Shoikhet's $L_\infty$-morphism $\mathcal{U}^{sh}$. 

%We now construct an explicit $L_\infty$-quasi-inverse $\mathcal U^{inv}$ to \eqref{eq:qiso2}. We keep the notations 
%introduced in the proof of Lemma \ref{lem:Poincare}. $\mathcal U_k^{inv}$ will be expressed via universal 
%expressions involving $\kappa$, $\iota$, $p$; and $[\cdot,\cdot]_{\rm SN}$. 

%[...]

%Therefore the first property is also obviously satisfied by $\mathcal U^{inv}$, as all the ingerdient 
%involved in its definition are defined for $\mathbb{R}^d_{\rm formal}$ as well. 

%The second and third properties are obviously satisfied by $\mathcal{U}^{cyc}$. Namely, 
%$\mathcal U^{inv}\circ\iota={\rm id}$, hence for any $\alpha_1,\dots,\alpha_k\in\mathfrak{g}$ one has 
%$$
%\mathcal{U}^{cyc}_k(\iota(\alpha_1),\dots\iota(\alpha_k))=\mathcal U^{sh}_k(\alpha_1,\dots,\alpha_k)\,.
%$$ 
%[...]

In the rest of the Section we follow closely \cite{dlg05} (see also \cite{dolgushev, dolgushev-2005}) which we adapt to the context of the cyclic formality. 
We assume the reader is familiar with the methods therein. 

\subsubsection{Recollection about Dolgushev-Fedosov resolutions}

We will consider differential forms with values in the following (graded) sheaves from \cite{dlg05, dolgushev, dolgushev-2005}: 
\begin{itemize}
\item the sheaf $\mathcal O$ of fiberwisely formal functions on $TM$; 
\item the sheaf $\mathcal T_{\rm poly}^{\bullet}$ of fiberwisely formal polyvector fields tangent to the fibers; 
\item the sheaf $\mathcal D_{\rm poly}^{\bullet}$ of fiberwisely formal polydifferential operators tangent to the fibers;
\item the sheaf $\mathcal A^{\bullet}$ of fiberwisely formal differential form tangent to the fibers. 
\end{itemize}
All these sheaves are acted on upon the (sheaf of) Lie algebra $\mathcal T:=\mathcal T_{\rm poly}^1$. 

\medskip

Now consider a torsion free connection $\nabla$ on $M$ with values in the tangent bundle and we let 
$\mathcal B$ be any of the four previously mentioned sheaves. 
Let us introduce local coordinates $(x^1,\dots,x^d)$ and write $(y^1,\dots,y^d)$, $y^i:={\rm d}x^i$, 
for the corresponding fiberwise coordinates on $TM$. 
Recall from \cite{dolgushev} that one can construct a differential $D_{\nabla}$ on 
$\Omega^\bullet(\mathcal B^\bullet)$ of the form 
\begin{equation}\label{eq:FeD}
D_{\nabla}={\rm d}x^i\frac{\partial}{\partial x^i}
+\big(\underbrace{A-{\rm d}x^i\Gamma_{ij}^ky^j\frac{\partial}{\partial y^k}
-{\rm d}x^i\frac{\partial}{\partial y^i}}_{=:Q}\big)\cdot\,,
\end{equation}
where $\Gamma_{ij}^k$ are Christoffel symbols of $\nabla$ and $A\in\Omega(M,\mathcal T)$, with the following properties:
\begin{itemize}
\item $D_\nabla$ respects all the algebraic structure on $\mathcal B$. E.g.~the (fiberwise) Schouten-Nijenhuis product 
on $\mathcal T_{\rm poly}^\bullet$, the (fiberwise) Gerstenhaber bracket and Hochschild differential on 
$\mathcal D_{\rm poly}^\bullet$, the natural (fiberwise) pairing between $\mathcal T_{\rm poly}$ and $\mathcal A$, the 
(fiberwise) de Rham differential on $\mathcal A$; 
\item The cohomology of $D_\nabla$ is concentrated in degree zero: 
\begin{itemize}
\item $H^\bullet\big(\Omega(M,\mathcal O),D_\nabla\big)=C^\infty(M)$ as algebras; 
\item $H^\bullet\big(\Omega(M,\mathcal T_{\rm poly}),D_\nabla\big)=T_{\rm poly}(M)$ as graded Lie algebras;
\item $H^\bullet\big(\Omega(M,\mathcal D_{\rm poly}),D_\nabla\big)=D_{\rm poly}(M)$ as DG Lie algebras;
\item $H^\bullet\big(\Omega(M,\mathcal A),D_\nabla\big)=\Omega(M)$ as DG algebras. 
\end{itemize}
\end{itemize}

Moreover, one can produce an explicit isomorphism 
$$
\lambda\,:\,{\rm B}^\bullet\,\longrightarrow\,\Omega^0(M,\mathcal B^\bullet)\cap\ker(D_\nabla)\,,
$$
of algebras (resp.~DG algebras, graded Lie algebras, DG Lie algebras), with ${\rm B}$ being $C^\infty(M)$ 
(resp.~$\Omega(M)$, $T_{\rm poly}M$, $D_{\rm  poly}M$) if $\mathcal B$ is $\mathcal O$ 
(resp.~$\mathcal A$, $\mathcal T_{\rm poly}$, $\mathcal D_{\rm poly}$). 

\medskip 

The resulting injective quasi-isomorphism (of DG algebras or DG Lie algebras) 
$$
{\rm B}^\bullet\,\longrightarrow\,\Omega^\bullet(M,\mathcal B^\bullet)
$$
is called the Dolgushev-Fedosov resolution of ${\rm B}^\bullet$, and we still denote it $\lambda$. 

\subsubsection*{Explicit formul\ae}

Let us now describe explicitly the constructon of $A$ and $\lambda$. 

We start by defining the linear operator $\kappa\,:\,\Omega(M,\mathcal B)\,\to\,\Omega(M,\mathcal B)$ being 
$$
\frac{1}{k+l}y^i\langle\frac{\partial}{\partial x^i},\cdot\rangle
$$
on $k$-forms with values in sections of $\mathcal B$ that are $l$-polynomial in the fibers if $k+l>0$ 
(and zero otherwise). It follows from a straightforward computation that $\kappa$ is a chain homotopy for 
$\delta:={\rm d}x^i\frac{\partial}{\partial y^i}\cdot$:
$$
\big(\delta\circ\kappa+\kappa\circ\delta\big)(s)=s-s_{|y^i={\rm d}x^i=0}
\qquad\big(\forall s\in\Omega(M,\mathcal B)\big)\,, 
$$
and $\kappa\circ\kappa=0$. 

Following \cite{dlg05} (see also \cite{dolgushev, dolgushev-2005}), one defines $A$ recursively as follows: 
\begin{equation}\label{eq:recu1}
A:=\kappa\left(-\frac12{\rm d}x^i\wedge{\rm d}x^j{R_{ij}}_k^ly^y\frac{\partial}{\partial y^l}
+{\rm d}x^i\left(\frac{\partial A}{\partial x^i}-[\Gamma_{ij}^ky^j\frac{\partial}{\partial y^k},A]^f\right)
+\frac12[A,A]^f\right)\,,
\end{equation}
where ${R_{ij}}_k^l$ is the curvature tensor of $\nabla$. This way, $A$ is such that $\kappa(A)=0$ and 
$A\equiv 0~\textrm{mod}~|y|^2$. Moreover, it is the unique such. 

\medskip

We define $\lambda$ similarly for ${\rm B}\neq D_{\rm poly}M$ (remark that in the cases under consideration 
there is a canonical identification between the space ${\rm B}$ and the space of sections of $\mathcal B$ that 
are constant in the fibers\footnote{It simply maps $\frac{\partial}{\partial x^i}$ to 
$\frac{\partial}{\partial y^i}$ and ${\rm d}x^i$ to ${\rm d}y^i$. The case of (poly)differential operators 
is more intricate. }). 
For any section $s_0$ of ${\rm B}$, $\lambda(s_0)$ is defined as the unique section $s$ of $\mathcal B$ 
defined by the following recursive relation: 
\begin{equation}\label{eq:recu2}
s:=s_0+\kappa\left({\rm d}x^i\left(\frac{\partial s}{\partial x^i}
-\Gamma_{ij}^ky^j\frac{\partial}{\partial y^k}\cdot s\right)+A\cdot s\right)\,.
\end{equation}

\subsection{Dolgushev-Fedosov resolution in the cyclic context}

Let $\omega$ be the volume form on $M$ and define $\widetilde\omega:=\lambda(\omega)$. We then have an operator 
$$
\widetilde{\dv}:={\dv}^f_{\widetilde{\omega}}\,:\,\mathcal T_{\rm poly}^\bullet\,
\longrightarrow\,\mathcal T_{\rm poly}^{\bullet-1}
$$
obtained as the composed map 
$$
\mathcal T_{\rm poly}^\bullet\,\overset{\langle\widetilde\omega,\cdot\rangle}{\longrightarrow}\,\mathcal A^{n-\bullet-1}\,
\,\overset{{\rm d}_{\rm dR}^f}{\longrightarrow}\,\mathcal A^{n-\bullet}\,
\overset{\langle\widetilde\omega,\cdot\rangle^{-1}}{\longrightarrow}\,\mathcal T_{\rm poly}^{\bullet-1}\,.
$$
\begin{lemma}
$\lambda\big(\dv_\omega(\alpha)\big)=\widetilde{\dv}\big(\lambda(\alpha)\big)$ for any 
$\alpha\in T_{\rm poly}M$. 
\end{lemma}
\begin{proof}
This is a direct consequence of the fact that $\lambda$ is a morphism of DG algebras (see above) 
and the obvious fact that $\lambda$ respects the pairing between polyvectors and forms. 
\end{proof}

From now and in the rest of the Section we assume that $\nabla$ preserves the volume form $\omega$ 
(i.e.~the covariant derivative of $\omega$ w.r.t.~$\nabla$ vanishes). Let us write $\omega$ in coordinates:  
$$
\omega=g(x^1,\dots,x^d){\rm d}x^1\wedge\cdots\wedge{\rm d}x^d\,.
$$
\begin{prop}\label{prop:divA}
Under the above assumption,  
$$
\widetilde\omega=g(x^1,\dots,x^d){\rm d}y^1\wedge\cdots\wedge{\rm d}y^d
$$
(i.e.~$\lambda(\omega)=\omega$), and thus $\widetilde{\dv}$ is the divergence operator defined 
by the standard volume form in the fibers. 

Moreover $\widetilde{\dv}(A)=0$, and thus $\widetilde{\dv}$ commutes with $D_\nabla$ 
(locally, $\widetilde{\dv}(Q)=0$). 
\end{prop}
\begin{proof}
The first statement of the proposition directly follows from the recursive definition \eqref{eq:recu2} 
for $\lambda(\omega)$. 

For the second statement we first observe that, since $\widetilde\omega$ is constant in the fibers, 
$\widetilde{\rm div}$ commutes with $\kappa$. Then as $\nabla$ preserves the volume form then one can see that 
$$
\widetilde{\rm div}\left({\rm d}x^i\wedge{\rm d}x^j{R_{ij}}_k^ly^y\frac{\partial}{\partial y^l}\right)=0\,.
$$
Hence it follows from the recursive definition \eqref{eq:recu1} of $A$ that $\widetilde{\rm div}(A)=0$. 

Finally, the last statement follows from  
$\widetilde{\dv}\left({\rm d}x^i\frac{\partial}{\partial y^i}\right)=0$ and the fact that 
$\nabla$ preserves the volume form. 
\end{proof}
Therefore the $u$-linear extension of $\lambda$ defines a quasi-isomorphism of DGLAs 
\begin{equation}\label{eq:reso1}
\big(T_{\rm poly}(M)[u], u \dv_\omega,[\cdot,\cdot]_{\rm S}\big)\,\longrightarrow\,
\Big(\Omega\big(M,\mathcal T_{\rm poly}[u]\big),D_\nabla+u \widetilde{\dv},[\cdot,\cdot]_{\rm S}^f\Big)\,.
\end{equation}
\bigskip

We now define $\mathcal{CC}:=(\widetilde{\mathcal{CC}})_{\mathcal T}$ as the (graded) sheaf of 
$\mathcal T$-coinvariants of 
$\widetilde{\mathcal{CC}}^\bullet:=\mathcal D_{\rm poly}^{\bullet+1}\otimes_{\mathcal O}\mathcal A^d$. 
One can prove along the same lines as in Section \ref{sec:hochcyc} that the composed map 
$$
\iota^f\,:\,\mathcal D_{\rm poly}^\bullet\,\hookrightarrow\,\mathcal D_{\rm poly}^{\bullet+1}\,
\overset{\cdot\otimes\widetilde\omega}{\longrightarrow}\,\widetilde{\mathcal{CC}}^\bullet\,
\twoheadrightarrow\,\mathcal{CC}^\bullet
$$
is an isomorphism of graded sheaves. 

Again, one has that $(\lambda\otimes\lambda)\circ\iota=\iota^f\circ\lambda$ 
(here we have used that $\lambda(1)=1$). Moreover, $\lambda$ also commutes with the cyclic shift operator $\sigma$: 
$\lambda\circ\sigma=\sigma^f\circ\lambda$, and it follows from Proposition \ref{prop:divA} that $D_\nabla$ preserves 
$(D_{\rm poly}(M))^\sigma$. Therefore, $\lambda$ restricts to a quasi-isomorphism of DGLAs 
\begin{equation}\label{eq:reso2}
\big((D_{\rm poly}(M))^{\sigma},{\rm d}_{\rm H},[\cdot,\cdot]_{\rm G}\big)\,\longrightarrow\,
\Big(\Omega\big(M,(\mathcal D_{\rm poly})_{\rm cyc}\big),D_\nabla+{\rm d}_{\rm H}^f,[\cdot,\cdot]_{\rm G}^f\Big)\,,
\end{equation}
where $(\mathcal D_{\rm poly})_{\rm cyc}$ is the graded subspace of $\mathcal D_{\rm poly}$ consisting of elements 
$P$ such that 
$$
\sigma^f\big(\iota^f(P)\big)=\iota^f(P)\,.
$$

\subsubsection{Proof of Theorem \ref{thm:cycform}}

Take $(U_\alpha,\rho_\alpha)$ a partition of unity with $U_\alpha$ beging open coordinate charts of $M$. 
Let us first apply fiberwisely the cyclic formality $L_\infty$-quasi-isomorphism over $U_\alpha$, 
and then consider its $\big(\Omega(U_\alpha),{\rm d}_{\rm dR}\big)$-linear extension: 
$$
\mathcal{U}^{cyc,\alpha}\,:\,
\left(\Omega\Big(U_\alpha,\mathcal T_{\rm poly}[u]\Big),{\rm d}x^i\frac{\partial}{\partial x^i}+u\widetilde{\dv},
[\cdot,\cdot]_{\rm SN}^f\right)
\,\longrightarrow\,
\left(\Omega\Big(U_\alpha,(\mathcal D_{\rm poly})_{\rm cyc}\Big),{\rm d}x^i\frac{\partial}{\partial x^i}+{\rm d}_{\rm H},
[\cdot,\cdot]_{\rm G}^f\right)\,.
$$
Let $Q_\alpha\in\Omega^1(U_\alpha,T_{\rm poly}^0)$ be the restriction of $Q$ to $U_\alpha$, expressed in coordinates. 
\begin{rem}
Let us remind to the reader that $Q$ is NOT a tensor (and {\it a fortiori} a vector field valued one-form), 
but that $Q_\alpha$ and $Q_\beta$ differ by a linear vector field (valued one-form) on $U_\alpha\cap U_\beta$. 
\end{rem}
It follows from \eqref{eq:FeD} and $D_\nabla\circ D_\nabla=0$ that $Q_\alpha$ is a Maurer-Cartan element: 
$$
{\rm d}x^i\frac{\partial Q_\alpha}{\partial x^i}+\frac12[Q_\alpha,Q_\alpha]_{\rm SN}^f=0\,.
$$
Since $Q_\alpha$ is a $\widetilde{\dv}$ergence free vector field (valued one-form) then 
$$
\widetilde{Q_\alpha}:=\sum_{n\geq1}\frac1{n!}\mathcal U_n^{cyc,\alpha}(Q_\alpha,\dots,Q_\alpha)=Q_\alpha\,,
$$
and thus the structure maps 
$$
(\mathcal U^{cyc,\alpha}_{Q_\alpha})_n(s_1,\dots,s_n):=
\sum_{k\geq0}\frac1{k!}\mathcal U^{cyc,\alpha}(\underbrace{Q_\alpha,\dots,Q_\alpha}_{k\textrm{ times}},s_1,\dots, s_n)
$$
define an $L_\infty$-quasi-isomorphism $\mathcal U^{cyc,\alpha}_{Q_\alpha}$ 
$$
\left(\Omega\Big(U_\alpha,\mathcal T_{\rm poly}[u]\Big),{\rm d}x^i\frac{\partial}{\partial x^i}+Q_\alpha\cdot
+u\widetilde{\dv},[\cdot,\cdot]_{\rm SN}^f\right)
\,\longrightarrow\,
\left(\Omega\Big(U_\alpha,(\mathcal D_{\rm poly})_{\rm cyc}\Big),{\rm d}x^i\frac{\partial}{\partial x^i}
+Q_\alpha\cdot+{\rm d}_{\rm H},
[\cdot,\cdot]_{\rm G}^f\right)\,.
$$
We then construct an $L_\infty$-quasi-isomorphism 
$\mathcal U^{cyc}_{\nabla}:=\sum_\alpha\rho_\alpha\mathcal U^{cyc,\alpha}_{Q_\alpha}$ 
$$
\left(\Omega\Big(M,\mathcal T_{\rm poly}[u]\Big),D_\nabla+u\widetilde{\dv},[\cdot,\cdot]_{\rm SN}^f\right)
\,\longrightarrow\,
\left(\Omega\Big(M,(\mathcal D_{\rm poly})_{\rm cyc}\Big),D_\nabla+{\rm d}_{\rm H},
[\cdot,\cdot]_{\rm G}^f\right)\,.
$$
{\it A priori}, $\mathcal U^{cyc}_{\nabla}$ is not well-defined. But since $Q_\alpha$ and $Q_\beta$ differ by a 
linear vector field (valued one-form) on intersections which is divergence free, then thanks to the last property 
of Proposition \ref{prop:properties} there is no ambiguity in the definition of $\mathcal U^{cyc}_{\nabla}$. 

\medskip

Finally, we obtained a chain of $L_\infty$-quasi-isomorphism 
$$
T_{\rm poly}(M)[u]\,\overset{\lambda}{\longrightarrow}\,\Omega\big(M,\mathcal T_{\rm poly}[u]\big)\,
\overset{\mathcal U^{cyc}_{\nabla}}{\longrightarrow}\,\Omega\big(M,(\mathcal D_{\rm poly})_{\rm cyc}\big)\,
\overset{\lambda}{\longleftarrow}(D_{\rm poly}(M))^{\sigma}\,.
$$
Hencefore we have proved Theorem \ref{thm:cycform}. \hfil \qed

\section{Application: Classification of closed star products}
\label{sec:applications}
This application of cyclic formality has already been proposed by Shoikhet in \cite{shoikhet}. We repeat it here briefly.
\begin{defi}
Let $M$ be a smooth manifold. We call an associative $\hbar$-linear product $\star$ on $C^\infty(M)[[\hbar]]$ a \emph{star product} if 
\[
f \star g = fg+\sum_{j\geq 1} \hbar^j m_j(f,g) 
\]
for bidifferential operators $m_j$ and all $f,g\in C^\infty(M)$. Two star products $\star$, $\star'$ are \emph{gauge equivalent} if there is a formal series of differential operators
\[
D = 1 + \sum_{j\geq 1} \hbar^j D_j
\]
such that $f\star g=D^{-1}(D(f)\star' D(g))$ for all $f,g\in C^\infty(M)$.
Let now $M$ be oriented with volume form $\omega$. The star product $\star$ is \emph{closed} if for any three compactly supported $f,g,h\in C_c^\infty(M)$
\[
\int_M f(g\star h) \omega = \int_M g(h\star f)\omega. 
\]
Two closed star products $\star$, $\star'$ are called \emph{gauge equivalent} if there is a $D$ as above such that $f\star g=D^{-1}(D(f)\star' D(g))$ and 
\[
\int_M fg \omega = \int_M D(f)D(g) \omega
\]
for all $f,g\in C_c^\infty(M)$.
\end{defi}

\begin{defi}
A formal series of degree $1$ elements of $T_{poly}^{\bullet+1}[u]$
\[
\pi = \sum_{j\geq 0} \hbar^j (\pi_j+u f_j)
\]
is called \emph{formal unimodular Poisson structure} if $u \dv \pi +\frac{1}{2} \co{\pi}{\pi}_{S}=0$. Two such $\pi$, $\pi'$ are \emph{gauge equivalent} if there is a formal series of vector fields 
\[
\xi = \sum_{j\geq 1} \hbar^j \xi_j
\]
such that 
\[
\pi' = e^{ad_\xi} \pi + u \frac{1-e^{ad_\xi}}{ad_\xi} \dv\xi
\]
where $ad_\xi (\cdot):= \co{\xi}{\cdot}_{S}$.
\end{defi}

The above notions of gauge equivalence obviously define equivalence relations on the sets of closed star products and formal unimodular Poisson structures. There is a map $\Phi$ between the two sets, mapping a formal unimodular Poisson structure $\pi$ to the (closed) star product $\star$ such that
\[
f\star g = fg + \sum_{j\geq 1} \frac{\hbar^j}{j!} \morphU_j(\pi\odot \cdots \odot \pi)(f,g) 
\]

\begin{thm}
The map $\Phi$ is a bijection between the set of gauge equivalence classes of formal unimodular Poisson structures and the set of gauge equivalence classes of closed star products. 
\end{thm}
\begin{proof}
Kontsevich has shown in \cite{kontsevich} that an $L_\infty$-quasi-isomorphism between dglas $\alg{g}\to \alg{g}'$ induces in the above manner a bijection between the gauge equivalence classes of solutions to the Maurer-Cartan equations in the pro-nilpotent Lie algebras $\hbar \alg{g}[[\hbar]]$ and $\hbar \alg{g}'[[\hbar]]$. Applying this to $\alg{g}=T_{poly}(M)[u]$ and $\alg{g}'=(D_{poly}(M))^\sigma$ yields the result.
\end{proof}

\begin{rem}
In case the top cohomology of $M$ vanishes, i.e., $H^d(M)=0$, one can see that any formal unimodular Poisson structure is gauge equivalent to a divergence free Poisson structure.
\end{rem}

\nocite{*}
\bibliographystyle{plain}
\bibliography{cycform} 

\end{document}